\documentclass[11pt]{amsart}
\newcommand{\remove}[1]{}

\numberwithin{equation}{section}
\usepackage{amssymb}

\usepackage{mathrsfs}

\usepackage{graphicx}
\usepackage{amsmath}
\usepackage{amsmath,amssymb}
\usepackage{yfonts}
\usepackage{mathtools}
\usepackage{mathabx}
\usepackage{enumitem}

 \DeclareMathAlphabet{\mathcalligra}{T1}{calligra}{m}{n}
 
 \DeclareMathOperator{\Char}{char}

  \DeclareMathOperator{\depth}{depth}

 \DeclareMathOperator{\grade}{gr}
 \DeclareMathOperator{\height}{height}

 \DeclareMathOperator{\orbit}{orbit}

 \DeclareMathOperator{\Span}{span}
 \DeclareMathOperator{\spec}{spec}

\newtheorem{Theorem}{Theorem}[section]
\newtheorem{Theorem--Definition}[Theorem]{Theorem--Definition}
\newtheorem{Corollary}[Theorem]{Corollary}
\newtheorem{Lemma}[Theorem]{Lemma}
\newtheorem{Lemma--Definition}[Theorem]{Lemma--Definition}

\newtheorem *{theorem}{Theorem}
\newtheorem *{theorema}{Theorem A}
\newtheorem *{theoremb}{Theorem B}
\newtheorem *{theoremba}{Theorem B*}
\newtheorem *{theoremc}{Theorem C}
\newtheorem *{theoremd}{Theorem D}
\newtheorem *{theoremda}{Theorem D*}
\newtheorem *{theoreme}{Theorem E}
\newtheorem *{theoremf}{Theorem F}
\newtheorem *{theoremg}{Theorem G}
\newtheorem *{theoremh}{Theorem H}
\newtheorem *{theoremi}{Theorem I}
\newtheorem *{theoremj}{Theorem J}
\newtheorem *{theoremk}{Theorem K}
\newtheorem *{theoreml}{Theorem L}
\newtheorem *{theoremm}{Theorem M}
\newtheorem *{question}{Question}
\newtheorem *{definition}{Definition}
\newtheorem *{remark}{Remark}
\newtheorem *{notation}{Notation}
\newtheorem *{note}{Note}
\newtheorem *{conjecture}{Conjecture}
\newtheorem *{conjecturec}{Conjecture C}
\newtheorem *{convention}{Convention}
\newtheorem *{warning}{Warning}

\newtheorem{Definition-Remark}[Theorem]{Definition-Remark}
\newtheorem{Proposition}[Theorem]{Proposition}
\newtheorem{Proposition--Definition}[Theorem]{Proposition--Definition}

\newtheorem{Note}[Theorem]{Note}

\theoremstyle{remark}

\theoremstyle{definition}
\newtheorem{Remark}[Theorem]{Remark}
\newtheorem{Example}[Theorem]{Example}
\newtheorem{Remark--Definition}[Theorem]{Remark--Definition}
\newtheorem{Definition}[Theorem]{Definition}

\begin{document}
\title{$p$-Groups and the Polynomial Ring of Invariants Question}
\author{Amiram Braun}
\vspace{.4in}

\date{\today\vspace{-0.5cm}}
\maketitle

\begin{center}
Dept. of Mathematics, University of Haifa, Haifa, Israel 31905

abraun@math.haifa.ac.il
\end{center}

\begin{abstract}
Let $G \subset GL(V)$ be a finite $p$-group, where $p=\Char F$ and $\dim _FV$ is finite.
Let $S(V)$ be the symmetric algebra of $V$, $S(V)^G$ the subring of $G$-invariants, and $V^*$ the dual space of $V$.
The following presents our solution to the above question.

\begin{theorema}
Suppose $\dim _FV=3$. Then $S(V)^G$ is a polynomial ring if and only if $G$ is generated by transvections.
\end{theorema}

\begin{theoremb}
Suppose $\dim _FV \ge 4$. Then $S(V)^G$ is a polynomial ring if and only if:
\begin{enumerate}
  \item $S(V)^{G_U}$ is a polynomial ring for each subspace $U \subset V^*$ with $\dim _F U=2$, where $G_U=\{g \in G | g(u)=u$, $\forall u \in U\}$, and
  \item $S(V)^G$ is Cohen-Macaulay.
\end{enumerate}
Alternatively, $(1)$ can be replaced by the equivalent condition:
\begin{enumerate}
 \setcounter{enumi}{2}
   \item $\dim ($ non-smooth locus of $S(V)^G)$ $\le 1$.
\end{enumerate}
\end{theoremb}
\end{abstract}

\section{\bf Introduction}\label{SEC-1}
Let $V$ be a finite dimensional vector space over a field $F$.
Let $G\subset GL(V)$ be a finite group and $S(V)$ the symmetric algebra of $V$.
$G$ acts as a finite group of automorphisms on $S(V)$.
Let $S(V)^G=\{x\in S(V)|g(x)=x,\forall g\in G\}$ be the ring of $G$-invariants.
The determination of the polynomial property of $S(V)^G$ is an old problem.
The basic result is the following.
We refer to \cite{BE}, \cite{NS} for an historical account.

\begin{theorem} (Shephard-Todd-Chevalley-Serre).
Suppose in addition that $(p,|G|)=1$, where $p=\Char F$. Then $S(V)^G$ is a polynomial ring if and only if $G$ is generated by pseudo-reflections.
\end{theorem}

Serre also proved the following implications, where no restrictions on $p$ are assumed.
\begin{theorem} \cite{Serre}, \cite[V, \S6, exercise 8]{BOU-AC}.
Suppose $S(V)^G$ is a polynomial ring.
Then $G$ is generated by pseudo-reflections and $S(V)^{G_U}$ is a polynomial ring for each subspace $U \subset V^*$, where $G_U=\{g\in G |g(u)=u, \forall u \in U\}$.
\end{theorem}

The second implication was firstly proved by R. Steinberg \cite{St} in case $F=\mathbb{C}$.
This theorem suggests a starting point in answering the polynomiality question in the modular case, that is if $p$ divides $|G|$.
The following result by Nakajima was the first positive result.

\begin{theorem} \cite{N2}.
Suppose $G\subset GL(V)$ is a finite $p$-group, where $\dim_FV=n$ is finite, and $F=\mathbb{F}_p$ is the prime field.
Then $S(V)^G$ is a polynomial ring if and only if:
\\There is a basis of $V$ $\{z_1,...,z_n\}$ with $g(z_1)=z_1$, $(g-1)(z_i)\in Fz_1+\cdots +Fz_{i-1}$, $i=2,...,n$, for each $g \in G$ such that:
$$|G-\orbit (z_1)|\cdots |G-\orbit (z_n)|=|G|.$$
\end{theorem}
Actually the "if" part always implies the polynomiality of $S(V)^G$ with no restriction on $F$.
The necessity of the condition is only granted when $F=\mathbb{F}_p$.
As such this theorem has no relation to Serre's result.

Extending and completing earlier results of Nakajima, Kemper-Malle established,using classification theorems for finite irreducible groups, the following result.

\begin{theorem} (\cite[Main Theorem]{KM} and the corrections in \cite[p. 119]{DK}).
Suppose $G \subset GL(V)$ is a finite group, where $V$ is a finite dimensional irreducible $G$-module.
Then $S(V)^G$ is a polynomial ring if and only if:
\begin{enumerate}
\item $S(V)^{G_U}$ is a polynomial ring for each subspace $U \subset V^*$ with $\dim_FU \ge 1$, and
\item $G$ is generated by pseudo-reflections.
\end{enumerate}
\end{theorem}

Actually one merely needs assuming $\dim _F U=1$ in item \textit{(1)}.
However, in spite of a reduction to an indecomposable $V$, Kemper-Malle refrained from stating the following conjecture.

\begin{conjecturec}
Suppose $G \subset GL(V)$ is a finite group, where $V$ is finite dimensional.
Then $S(V)^G$ is a polynomial ring if and only if:\\
$\indent (C1)$ $S(V)^{G_U}$ is a polynomial ring for each subspace $U \subset V^*$ with $\dim_FU= 1$, and
$\indent (C2)$ $G$ is generated by pseudo-reflections.
\end{conjecturec}

Our result, in case of $\dim_FV\ge 4$, initiated from the trivial observation that Conjecture C holds trivially for a $p$-group $G$.
Indeed the action of $G$ on $V^*$ can always be simultaneously triangulized, in particular there exists $f \in V^*$, $f \neq 0$ with $g(f)=f$, for each $g \in G$.
Therefore $G_U=G$, where $U=Ff$ and Conjecture C trivially holds.
This suggests, in the spirit of Conjecture C, that for a finite $p$-group, there may exist a weaker assumption on subspaces of $V^*$, e.g. on  the set $\{U|U \subset V^*$, $\dim_F U=2\}$.
Indeed our main result is of this nature:

\begin{theoremb}
Suppose $G \subset GL(V)$ is a finite $p$-group, where $p=\Char F$, $\dim_FV \ge 4$. Then $S(V)^G$ is a polynomial ring if and only if:
\begin{enumerate}
  \item $S(V)^{G_U}$ is a polynomial ring, for each subspace $U \subset V^*$, with $\dim _F U=2$, where $G_U=\{g \in G | g(u)=u$, $\forall u \in U\}$, and
  \item $S(V)^G$ is Cohen-Macaulay.
\end{enumerate}
\end{theoremb}
Moreover \textit{(1)} is equivalent to:
\\ \indent \textit{(3)} $\dim ($\textit{non-smooth locus of} $S(V)^G)$ $\le 1$.

So Theorem B has a singularity theory flavour. Note that the assumption $\dim_FV \ge 4$ is essential, as every finite $p$-subgroup of $GL_3(F)$ satisfies conditions B(1) and B(2).

One may be puzzled by the omission of the condition "$G$ is generated by pseudo-reflections", in Theorem B.
The simple answer is that conditions B(1) and B(2) actually imply it.
This is a consequence of a result of Kemper \cite[Corollary 3.7]{K-CM}.

Another seemingly puzzling question is the appearance of the Cohen-Macaulay assumption in Theorem B(2), as it neither appears in Kemper-Malle \cite{KM} result nor in Conjecture C.
In fact the Cohen-Macaulay property is a consequence of condition (C1).
Indeed condition (C1) implies $\dim $(singular locus of $S(V)^G)=0$, and hence $\dim $(non-Cohen-Macaulay locus of $S(V)^G)=0$. Therefore by Kemper \cite[Thm. 3.1]{K-Loci} $S(V)^G$ is Cohen-Macaulay.

In case of $\dim_FV=3$ we have the following.

\begin{theorema}
Suppose $G \subset GL(V)$ is a finite $p$-group, where  $\dim _FV=3$ and $p=\Char F$.
Then $S(V)^G$ is a polynomial ring if and only if $G$ is generated by transvections.
\end{theorema}

The proof of Theorem A is computational.
By a standard result $G$ has the following presentations of matrices $G \subset
\begin{bsmallmatrix}
    1 & F & F \\
    0 & 1 & F \\
    0 & 0 & 1
  \end{bsmallmatrix}$,
namely $V$ has a basis $\{v_1,v_2,v_3\}$ with $g(v_1)=v_1$, $(g-1)(v_2)\in Fv_1$, $(g-1)(v_3)\in Fv_1+Fv_2$, for each $g\in G$.
The simpler cases $G \subset
\begin{bsmallmatrix}
    1 & F & F \\
    0 & 1 & 0 \\
    0 & 0 & 1
  \end{bsmallmatrix}$
and $G \subset
\begin{bsmallmatrix}
    1 & 0 & F \\
    0 & 1 & F \\
    0 & 0 & 1
  \end{bsmallmatrix}$
are handled by earlier results of Nakajima.
The remaining case is more involved and the end result $S(V)^G=F[v_1,a,b]$ where $a \in F[v_1,v_2]^G$, is crucial in the proof of Theorem B.

The following is a sketch of the proof of Theorem B, if $F$ is perfect (Theorem \ref{T2-12}).
It is done by a sequence of hypersurface sections reducing the problem, after $\dim_FV-3$ steps to a $3$-dimensional non-linear problem.
Theorem A cannot be used directly, but provides us with three elements which are a system of parameters (s.o.p.).
The homogenous lifts in $S(V)^G$ of these three elements together with the homogenous lifts in $S(V)^G$ of the elements involved in the definition of each hypersurface section, form a homogenous s.o.p. with their degree's product being equal to $|G|$.
This implies that $S(V)^G$ is generated by these $n$ lifts and consequently is a polynomial ring.
The general case requires some further considerations about the non-smooth locus.
A simplified proof in case $dim_FV=4$ is presented in the appendix.

We end with two $\dim_FV=4$ examples.
In the first one $|G|=p^3$, $G$ is elementary abelian and $S(V)^G$ is not a polynomial ring.
However $S(V)^G$ is Cohen-Macaulay, $G$ is generated by transvections, and $G=G_U$, $U \subset V^*$, with $\dim_F U=2$.
In the second example, we have, for the same group $G$, that $S(V^*)^G$ is a polynomial ring.

\begin{convention}
Matrices represent actions on the right.
\end{convention}

\begin{warning}
Since $g \in GL(V)$ is a pseudo-reflection on $V$ if and only if $g$ is a pseudo-reflection on $V^*$, we shall omit the reference to either $V$ or $V^*$.
Note however that the polynomiality of $S(V)^G$ or of $S(V^*)^G$ are not equivalent properties, e.g.  \cite[p. 120]{DK} or Examples \ref{EX2-1},  \ref{EX2-2}.
\end{warning}

\section{\bf Theorem B}\label{SEC-2}

\begin{notation}
Let $\mathfrak{m}$ be a prime ideal in $S(V)$.
Recall that $D_G(\mathfrak{m}):=\{g\in G|g (\mathfrak{m})=\mathfrak{m}\}$ is the decomposition subgroup of $\mathfrak{m}$ in $G$,
$I_G(\mathfrak{m}):=\{g \in G|(g-1)(S(V))\subseteq \mathfrak{m} \}$, the inertia subgroup of $\mathfrak{m}$ in $G$,
$W:=\mathfrak{m} \cap V$ and $U:=W^{\perp}=\{f \in V^*|f(W)=0\}$,
\end{notation}

\begin{Lemma}\label {L2-1}
Let $G \subset GL(V)$ be a finite group where $\dim _FV$  is finite. With the above notation we have:
$G_U=I_G(\mathfrak{m})$.
\end{Lemma}

\textbf{Proof:}
This follows from the observation that $G_U=\{g \in G|(g-1)(V) \subseteq W\}$.$\Box$

For a local commutative Noetherian ring $(A,\mathfrak{m})$, $\hat{A}$ denotes the completion of $A$ at $\mathfrak{m}$.

\begin{Theorem}\label {T2-2} (\cite[Expos\'{e} V, Prop. 2.2, p. 111]{G}).
With the notations above and setting $\mathfrak{n}:=\mathfrak{m} \cap S(V)^{I_G(\mathfrak{m})}$, $\mathfrak{p}:=\mathfrak{m}\cap S(V)^G$, $\mathfrak{q}:= \mathfrak{m} \cap S(V)^{D_G(\mathfrak{m})}$, we have:
\begin{enumerate}
\item $\widehat{S(V)_{\mathfrak{p}}^G}\cong \widehat{S(V)_{\mathfrak{q}}^{D_G(\mathfrak{m})}}$;
\item $S(V)_{\mathfrak{n}}^{I_G(\mathfrak{m})}$ is $\acute{e}$tale over  $S(V)_{\mathfrak{q}} ^{D_G(\mathfrak{m})}$.
\end{enumerate}
Consequently by \cite[the corollary in p. 179 and Theorem 23.7]{Matsumura} $S(V)^G_{\mathfrak{p}}$ is regular if and only if $S(V)_\mathfrak{n}^{I_G({\mathfrak{m}})}$ is regular.
\end{Theorem}

\begin{Remark}\label {R2-3}\
\begin{enumerate}
  \item There is a proof of Theorem \ref{T2-2}(1), for geometric points, in \cite[Prop. 1.1]{K-Loci}. However, it relies on \cite[Lemma 1.2]{K-Loci} which is false in the stated generality. This can be easily corrected in our setting, as was also observed by David Mundelius.
  \item Theorem \ref{T2-2} is proved in \cite{G} for any Noetherian ring $A$ with a finite subgroup $G \subset Aut(A)$.
\end{enumerate}
\end{Remark}

We have the following consequences, most of which also appear in \cite{K-Loci}. A further complement will be proved in Lemma \ref{L2-21}.

\begin{Lemma}\label {L2-4}
Let $G \subset GL(V)$ be a finite group and $\dim_FV=n$ is finite.
Consider the following conditions:
\begin{enumerate}
  \item $S(V)^G$ satisfies Serre's $R_m$ condition \cite[\S23, p. 183]{Matsumura};

  \item $S(V)^{G_U}$ is a polynomial ring for each subspace $U \subset V^*$ with $\dim_FU = n-m$;

  \item $S(V)^{G_U}$ is a polynomial ring for each subspace $U \subset V^*$ with $\dim_FU \ge n-m$.
\end{enumerate}
Then we have the following implications: $(2) \Leftrightarrow (3)$, $(3) \Rightarrow (1)$. Moreover if $F=F^p$ then $(1) \Rightarrow (3)$.
\end{Lemma}

\textbf{Proof:}
Let $\mathfrak{p}$ be a prime ideal in $S(V)^G$ with $\height(\mathfrak{p})=m$.
Let $\mathfrak{m}$ be a prime ideal in $S(V)$ with $\mathfrak{m} \cap S(V)^G=\mathfrak{p}$.
By "going up" we have $\height(\mathfrak{m})=\height(\mathfrak{p})=m$.
Let $W=V\cap \mathfrak{m}$.
Then $WS(V)\subseteq \mathfrak{m}$ is a prime ideal in $S(V)$ and $\dim _FW=\height(WS(V)) \le \height(\mathfrak{m})=m$.
Therefore $\dim _FU = \dim_F V- \dim _F W \ge n-m$, where $U:=W^\perp =\{f \in V^*|f(W)=0\}$.

We now show the implication $(3) \Rightarrow (1)$. By assumption $S(V)^{G_U}$ is a polynomial ring. Therefore by Lemma \ref{L2-1}, $S(V)^{I_G(\mathfrak{m})}$ is a polynomial ring and consequently $S(V)_{\mathfrak{n}}^{I_G(\mathfrak{m})}$ is regular where $\mathfrak{n}:=\mathfrak{m} \cap I_G(\mathfrak{m})$. So by Theorem \ref{T2-2} $S(V)_{\mathfrak{p}}^G$ is regular.

We show that $(2) \Leftrightarrow (3)$. It is clear that $(3) \Rightarrow (2)$. Assume $(2)$. Let $U' \subseteq V^*$ be a subspace with $\dim_FU' \ge n-m$. Let $U \subseteq U'$ be a subspace with $\dim_FU=n-m$. Therefore by assumption $S(V)^H$ is a polynomial ring, where $H:=G_U$. Consequently by Serre \cite[V, $\S$6, Ex.8]{BOU} $S(V)^{H_{U'}}$ is a polynomial ring. Let $g \in G_{U'}$. Then $g \in G_U=H$. Therefore $G_{U'}=H_{U'}$. Hence $S(V)^{G_{U'}}$ is a polynomial ring, as needed.

Assume now that $F=F^p$ and $(1)$ holds. Let $U \subseteq V^*$ be a subspace with $\dim_FU\ge n-m$. Set $W:=U^{\perp}=\{v \in V|f(v)=0, \forall f \in U\}$.
Then $\dim_FW=\dim_FV-\dim_FU \le n-(n-m)=m$.
Hence $WS(V)$ is a prime ideal with $\height(WS(V))=\dim_FW \le m$. Set $\mathfrak{p}=WS(V)\cap S(V)^G$. Then by "going down" $\height(\mathfrak{p})\le m$. Consequently by $(1)$ $S(V))_{\mathfrak{p}}^G$ is regular. Therefore by \cite[Prop. 1.9, Example 1.10 (where $F=F^p$ is mistakenly missing)]{K-Loci}, $S(V)^{I_G(WS(V))} =S(V)^{G_U}$ is a polynomial ring.$\Box$

We shall need the following observation.

\begin{Lemma}\label {L2-5}
Let $G \subset GL(V)$ be a finite $p$-group where $n:=\dim_ FV $ is finite and $p=\Char F$.
Suppose:
\begin{enumerate}
  \item $S(V)^G$ is Cohen-Macaulay;

  \item $S(V)^{G_U}$ is a polynomial ring for each subspace $U \subset V^*$ with $\dim_FU =n-2$.
\end{enumerate}
Then $G$ is generated by transvections.
\end{Lemma}

\textbf{Proof:} By $(1)$ and \cite[Cor. 3.7]{K-CM} we have that $G$ is generated by bi-reflections on $V^*$, $G=<\sigma _1,...,\sigma_r>$.
Since $\ker(\sigma _i -1):=U_i$, $U_i \subset V^*$ with $\dim_F U_i=n-2$, $i=1,...,r$, it follows by $(2)$ that $\sigma _i \in G_{U_i}$ where $S(V)^{G_{U_i}}$ is a polynomial ring. Therefore $G_{U_i}$ is generated by transvections for $i=1,...,r$, implying that the combined set of transvection generators of $G_{U_i}$, $i=1,...,r$, will generate $G$.$\Box$

The following is easy to verify.

\begin{Lemma}\label {L2-6}
Let $\sigma$ be a transvection on $V$ (equivalently on $V^*$) and $v \in V^{\sigma}$.
Then $\bar{\sigma}$, the induced map on $\bar{V}=V/Fv$ is a transvection (or the identity) on $\bar{V}$ (equivalently on $\bar{V}^*$).
\end{Lemma}

\textbf{Proof:} Clearly $\widebar{\ker (\sigma -1)} \subseteq \ker(\bar{\sigma} -1)$, so the result follows since $\dim _F\widebar{\ker(\sigma -1)}=\dim_F \widebar{V} -1$.$\Box$

The next result is a slight extension of Nakajima's \cite[Prop. 3.5]{N1}, \cite[Thm. 6.2.7]{NS} where the restriction $F=\mathbb{F}_q$ appears, but it is not used.

Let $A=\oplus_{n\ge0}A_n$, $A_0=F$, be a graded ring. We denote by $A_{+}:=\oplus_{n>0}A_n$, the irrelevant maximal ideal of $A$.

\begin{Proposition} \label{P2-6.5}
Let $\{g_1,...,g_k\} \subset GL(V)$, where $\dim _FV$ is finite and $\Char F >0$.
Let $\{v_1,...,v_r\}$ be a $F$-linearly independent subset of $V$ and $x\in S(V)-F[v_1,...,v_r]$, a homogenous element with $\deg x=p^e$.
Suppose $(g_j-1)(x) \in F[v_1,...,v_r]_{+}$, $g_j(v_i)=v_i$, $j=1,...,k$, $i=1,...,r$.
Then:
\begin{enumerate}
  \item $<g_1,...,g_k>$ induces an elementary abelian finite $p$-group of graded automorphisms on the polynomial ring $F[v_1,...,v_r,x]$.
  \item $F[v_1,...,v_r,x]^{<g_1,...,g_k>}=F[v_1,...,v_r,y]$ is a polynomial ring, $\deg y = p^t= \\p^e \cdot |<g_1,...,g_k>|$, and $y$ is a monic  $p$-polynomial in $x$ with coefficients in $F[v_1,...,v_r]_{+}$.
\end{enumerate}
\end{Proposition}

\textbf{Proof:}
Item (1) clearly holds.
We shall establish item (2) by induction on $k$.
Set $y_1:=x^p-[(g_1-1)(x)]^{p-1}x$ so $\deg y_1=p^{e+1}$, since $(g_1 - 1)(x)\in F[v_1,...,v_r]_{+}$, it follows that $(g_1-1)(y_1)=0$, as well as
$(g_j - 1)(y_1)=[(g_j-1)(x)]^p - [(g_1-1)(x)]^{p-1}(g_j-1)(x) \in F[v_1,...,v_r]_{+}$.
We claim that $F[v_1,...,v_r,x]^{<g_1>}=F[v_1,...,v_r,y_1]$.

By Galois theory and since $Q(F[v_1,...,v_r,x])^{<g_1>}=Q(F[v_1,...,v_r,x]^{<g_1>})$, we have $[Q(F[v_1,...,v_r,x]):Q(F[v_1,..,v_r,x]^{<g_1>})]=p$.
Since $\{1,x,...,x^{p-1}\}$ generate $F[v_1,...,v_r,x]$ as a $F[v_1,...,v_r,y_1]$-module, it follows that $\{1,x,...,x^{p-1}\}$ is a free basis of $F[v_1,...,v_r,x]$ over $F[v_1,...,v_r,y_1]$.

Let $z \in F[v_1,...,v_r,x]^{<g_1>}$. Then $z=c_1+c_2x+\cdots +c_{p-1}x^{p-1}$, $c_i\in F[v_1,...,v_r,y_1]$, $i=1,,,.p-1$. Let $s=\max \{i | c_i \neq 0 \}$. Suppose $s \ge 2$.

Recall that if $\phi  $ is an automorphism of a commutative ring $A$ and $\delta := \phi - 1$, then for $a \in A$ and $i \ge 1$ we have:
$$\delta (a^i)=i a^{i-1}\delta (a) + \sum_{j=0}^{i-2}\binom{i}{j}a^j \delta (a)^{i-j}.$$
Taking $\phi =g_1$, we get that $0 = \delta (z) =s c_s x^{s-1} \delta (x) + b$, where $b \in \sum_{j=0}^{s-2}F[v_1,...,v_r,y_1]x^j$.
Since $\delta (x) \neq 0$ (otherwise $g_1=1$), it follows from the freeness   of $\{1,x,...,x^{p-1}\}$ and  $
\delta (x)\in F[v_1,...,v_r]_{+}$  that $c_s=0$.
Hence $z\in F[v_1,...,v_r,y_1]$ and the case $k=1$ is established.
Observe that $K.\dim F[v_1,...,v_{r},y_1]=K.\dim F[v_1,...,v_r,x]^{<g_1>} =K.\dim F[v_1,...,v_r,x] = r+1$, implying that $F[v_1,...,v_r,y_1]$ is a polynomial ring in $r+1$ variables.

Since $(g_j-1)(y_1) \in F[v_1,...,v_r]_{+}$, $j=2,...,k$, as seen above, the general case follows now by induction, having by (1)
the $<g_2,...g_k>$-stability of $F[v_1,...,v_r,v_x]^{<g_1>}$ and the following:
$$F[v_1,...,v_r,x]^{<g_1,...,g_k>} =(F[v_1,...,v_r,x]^{<g_1>})^{<g_2,...,g_k>}
=F[v_1,...,v_r,y_1]^{<g_2,...,g_k>}.\Box$$

\begin{Remark} \label{R2-6.8}
\cite[Prop. 3.5]{N1} is obtained by taking $\{v_1,...,v_{n-1},v_n\}$ a basis of $V$, $x=v_n$ and $r=n-1$.
\end{Remark}

The next result extends to a non-linear setting Nakajima's \cite[Prop. 3.4, where $F=\mathbb{F}_q$ is unnecessarily assumed]{N1}. The proof is markedly different.

\begin{Proposition} \label{P2-7}
Let $\{g_1,...,g_k\} \subset GL(V)$, where $\dim_FV=n$ and $\Char F=p>0$. Let $\{v,m_2,...,m_n\} \subset S(V)$ satisfy:
\begin{enumerate}
\item $\{v,m_2,...,m_n\}$ are algebraically independent over $F$;
\item $v \in V$, $g_j(v)=v$, $j=1,...,k$;
\item $m_i$ is homogenous, $\deg m_i =p^{e_i}$ and $(g_j-1)(m_i) \in Fv^{p^{e_i}}$, $j=1,...,k$, $i=2,...,n$.
Then:
     \\(a) $<g_1,...,g_k>$ induces an elementary abelian finite $p$-group of graded automorphisms on the polynomial ring $F[v,m_2,...,m_n]$ and,
     \\(b) $F[v,m_2,...,m_n]^{<g_1,...g_k>}=F[v,y_2,...,y_n]$, is a polynomial ring where $y_i$ is homogenous, $i=2,...,n$.
\end{enumerate}
Moreover suppose:
\begin{enumerate}
  \setcounter{enumi}{3}
\item $F=F^p$  and,
\item $m_i-v_i^{p^{e_i}} \in vS(V)$ where $v_i\in V$, $i=2,...,n$.
\\Then:
\\(c) $y_i-u_i^{p^{f_i}} \in vS(V)$, where $u_i \in V$, $\deg y_i=p^{f_i}$, $i=2,...,n$
\end{enumerate}
\end{Proposition}

\textbf{Proof:} $<g_1,...g_k>$ extends to a group of graded automorphisms on $S(V)$ in the obvious way.
We have by (3) that $(g_j-1)(m_i)=a_{ji}v^{p^{e_i}}$, $a_{ji} \in F$, $i=2,...,n$, $j=1,...,k$.
This implies that $(g_j-1)^2(m_i)=(g_j-1)(a_{ji}v^{p^{e_i}})=0$, and $(g_j-1)^p(m_i)=0$.
Hence $g_{j|F[v,m_{1},...,m_{n}]}^p=1$ for $j=1,...,k$.
Also $(g_jg_l)(m_i) = g_j(m_i+a_{li}v^{p^{e_i}}) = m_i+a_{ji}v^{p^{e_i}}+a_{li}v^{p^{e_i}} = (g_lg_j)(m_i).$
Hence $<g_1,...,g_k>$ induces an elementary abelian finite $p$-group of graded automorphism on $F[v,m_{1},...,m_{n}]$, and (a) is verified.

Let $m_r$ be such that $g_1(m_r)\neq m_r$.
Equivalently $a_{1r} \neq 0$.
We have by Proposition \ref{P2-6.5} that
$F[v,m_r]^{<g_1>}=F[v,y]$, where $y:=m_r^p-(a_{1r}v^{p^{e_r}})^{p-1}m_r$.

We shall now handle the case $k=1$. We may assume that $e_2 \le e_3 \le \cdots \le e_n$.
Suppose $r$ is the $1^{st}$ index with $(g_1-1)(m_r)\neq 0$, that is $a_{1r} \neq 0$.
For $j\neq r$, $j \in \{2,...,n\}$ set $m'_j:=m_j - \frac{a_{1j}}{a_{1r}}v^{p^{e_j-e_r}}m_r \in F[v,m_2,...,m_r]$.
Clearly $m_j'=m_j$ if $j < r$.
Hence $(g_1-1)(m_j')=0$, $\forall j \neq r$.
Clearly $F[v,m_2,...,m_r]= F[v,m_r,m_j',j \neq r]$.
Consequently $\{v,m_r,m_2',...,m_{r-1}',m_{r+1}',...,m_n\}$ are algebraically independent over $F$.
\\Therefore $F[v,m_2,...,m_n]^{<g_1>} =F[v,m_r, m_j', j\neq r]^{<g_1>} =F[v,m_r]^{<g_1>}[m_j', j\neq r]= F[v,y,m_2',...,m'_{r-1},m'_{r+1},...,m'_n]$.
Since $n=K.\dim F[v,m_2,...,m_n] = K.\dim F[v,m_2,...,m_n]^{<g_1>}$, it follows that the latter is a polynomial ring, and the case $k=1$ is established.

We shall now settle the general case by induction on $k$.
We have by (a) that \\
$F[v,m_2,...,m_n]^{<g_1>}$ is $<g_2,...,g_k>$-stable, so it follows that:
\begin{equation*} 
\begin{split}
F[v,m_2,...,m_n]^{<g_1,...,g_k>} & = (F[v,m_2,...,m_n]^{<g_1>})^{<g_2,...,g_k>} \\
 & = F[v,y,m_2',....,m'_{r-1},m'_{r+1},...,m_n]^{<g_2,...g_k>}.
\end{split}
\end{equation*}
It is also evident from $m_j'=m_j-\frac{a_{1j}}{a_{1r}}v^{p^{e_j-e_r}}m_r$, that $\{g_2,...,g_k\}$ satisfy $(1)$, $(2)$ ,$(3)$ on $\{v,y,m_2',...,m'_{r-1},m'_{r+1},...,m'_n\}$, so by induction $F[v,y,m_j',j \neq r]^{<g_2,...,g_n>}$ is a polynomial ring.

We now need to verify item $(c)$. For this we need to show the existence of $\{u,u_j,j\neq r\}\subset V$ with $y-u^{\deg y} \in vS(V)$, $m'_j-u_j^{\deg m_j'} \in vS(V)$, $j\neq r$. Now if $e_j > e_r$ take $u_j:=v_j$ and $m_j'-u_j^{p^{e_j}}=(m_j-v_j^{p^{e_j}})-\frac{a_{1j}}{a_{1r}}v^{p^{e_j-e_r}}m_r \in vS(V)$.

If $e_j=e_r$, then $m_j'=m_j-\frac{a_{1j}}{a_{1r}}m_r$ is of degree $p^{e_r}$, so set $u_j:=v_j-(\frac{a_{1j}}{a_{1r}})^{\frac{1}{p^{e^r}}}v_r$.
Then $m_j'-u_j^{p^{e_r}}=(m_j - v_j^{p^{e_r}})-\frac{a_{1j}}{a_{1r}}(m_r-v_r^{p^{e_r}}) \in vS(V)$.
Now $u:=v-a_{1r}^{\frac{p-1}{p^{e_r+1}}}v \in V$ satisfies:
$y-u^{p^{e_r+1}}=[m_r^p-(a_{1r}v^{p^{e_r}})^{p-1}m_r] - [v^{p^{e_r+1}}-a_{1r}^{p-1}v^{p^{e_r+1}}]=(m_r-v^{p^{e_r}})^p - (a_{1r}v^{p^{e_r}})^{p-1}(m_r-v^{p^{e_r}}) \in vS(V)$.
This settles $(c)$ for $k=1$. The rest follows by induction as above.$\Box$

The next result is also an extension of Nakajima's \cite[Prop. 3.4]{N1}.

\begin{Corollary} \label {C2-8}
  Let $<g_1,...,g_k> \subset GL(V)$, where $\dim _FV=n$, $\Char F=p>0$ and $F=F^p$.
  Let $\{v,v_2,...,v_n\}$ be a basis of $V$ and assume that $g_j(v)=v$, $(g_j-1)(v_i) \in Fv$, $j=1,...,k$, $i=2,...,n$.
  Let $0 \le e_1 \le e_2 \le \cdots \le e_k$ be an increasing sequence of integers. Then:
  \begin{enumerate}
    \item $<g_1,...,g_k> \subseteq Aut_{\grade} (F[v^{p^{e_1}},v_2^{p^{e_2}},...,v_n^{p^{e_n}}]) \cap Aut_{\grade} (F[v,v_2^{p^{e_2-e_1}},...,v_n^{p^{e_n-e_1}}]) $;
    \item $F[v^{p^{e_1}},v_2^{p^{e_2}},...,v_n^{p^{e_n}}]^{<g_1,...,g_k>}=F[v^{p^{e_1}},z_2,...,z_n]$, and is a polynomial ring in $n$-variables, where $z_i$ is homogenous of degree $p^{h_i}$, $i=2,...,n$;
    \item $z_i-u_i^{p^{h_i}} \in v^{p^{e_1}}S(V)$, $i=2,...,n$ and $\{v,u_2,...,u_n\}$ is a basis of $V$;
    \item $e_1 \le h_2 \le \cdots \le h_n$.
  \end{enumerate}
\end{Corollary}

\textbf{Proof:}
$<g_1,...,g_k> \subset AutS(V)$ and the above assumptions grant that $F[v,v_2^{p^{e_2-e_1}},...,v_n^{p^{e_n-e_1}}]$ and $F[v^{p^{e_1}},v_2^{p^{e_2}},...,v_n^{p^{e_n}}]$ are $<g_1,...,g_k>$-stable.
Hence $(1)$ holds.

We next observe that the map $\phi (f) = f^{p^{e_1}}$, $\phi : F[v,v_2^{p^{e_2-e_1}},...,v_n^{p^{e_n-e_1}}] \rightarrow F[v^{p^{e_1}},v_2^{p^{e_2}},...,v_n^{p^{e_n}}]$, is a ring isomorphism (but not an algebra homomorphism since $F=F^p$ is being used).
It is also easily checked that $g_i \phi = \phi g_i$, $i=1,...,k$.
Consequently the restriction of $\phi $ induces the following isomorphism:
\begin{equation} \label{eq2:2}
\phi : F[v,v_2^{p^{e_2 - e_1}},...,v_n^{p^{e_n-e_1}}]^{<g_1,...,g_k>} \cong F[v^{p^{e_1}},v_2^{p^{e_2}},...,v_n^{p^{e_n}}]^{<g_1,...,g_k>}.
\end{equation}
Now by Proposition \ref{P2-7}, taking $m_i=v_i^{p^{e_i-e_1}}$, $i=2,...,n$, we have that
\\$F[v,v_2^{p^{e_2 - e_1}},...,v_n^{p^{e_n-e_1}}]^{<g_1,...,g_k>} = F[v,y_2,...,y_n]$ is a polynomial ring, where $y_i$ is homogenous of degree $p^{f_i}$, and $y_i-u_i^{p^{f_i}} \in vS(V)$, $i=2,...,n$.
Consequently by applying $\phi $ and using isomorphism (\ref{eq2:2}) we have that $F[v^{p^{e_1}},v_2^{p^{e_2}},...,v_n^{p^{e_n}}]^{<g_1,...,g_k>}=F[v^{p^{e_1}},y_2^{p^{e_1}},...,y_n^{p^{e_1}}]$,
so (2), (4) and part of (3) follow by setting $z_i:=y_i^{p^{e_1}}$, $z_i-u_i^{p^{f_i+e_1}} \in v^{p^{e_1}}S(V)$ and $h_i:=f_i+e_1$, $i=2,...,n$.

We shall now establish the basis property of $\{v,u_2,...,u_n\}$.
$v^{p^{e_1}}$ is a prime element in $F[v^{p^{e_1}},v_2^{p^{e_2}},...,v_n^{p^{e_n}}]$.
Since $S(V)=F[v,v_2,...,v_n]$ is finite over the former, it follows that $vS(V) \cap F[v^{p^{e_1}},v_2^{p^{e_2}},...,v_n^{p^{e_n}}]= v^{p^{e_1}}F[v^{p^{e_1}},v_2^{p^{e_2}},...,v_n^{p^{e_n}}]$.
Since $F[v^{p^{e_1}},v_2^{p^{e_2}},...,v_n^{p^{e_n}}]^{<g_1,...,g_k>}=F[v^{p^{e_1}},z_2,...,z_n]$ it follows that $v^{p^{e_1}}F[v^{p^{e_1}},v_2^{p^{e_2}},...,v_n^{p^{e_n}}]\cap F[v^{p^{e_1}},z_2,...,z_n] = v^{p^{e_1}}F[v^{p^{e_1}},z_2,...,z_n]$.
Consequently $vS(V) \cap F[v^{p^{e_1}},z_2,...,z_n]=v^{p^{e_1}}F[v^{p^{e_1}},z_2,...,z_n]:=(v^{p^{e_1}})$.
\\Note that $F[v^{p^{e_1}},z_2,...,z_n]/(v^{p^{e_1}})$ is a polynomial ring in $n-1$ variables.
Therefore:
\begin{equation*} 
\begin{split}
S(\bar{V})  = S(V)/vS(V) & \supset F[v^{p^{e_1}},z_2,...,z_n]/(v^{p^{e_1}}) = \\
& = F[\bar{z_2},...,\bar{z_n}]
 = F[\bar{u_2}^{p^{h_2}},...,\bar{u_n}^{p^{h_2}}],
\end{split}
\end{equation*}
where $\bar{V}:=V/Fv$ and $\bar{z_i}$ (respectively $\bar{u_i}$) is the image of $z_i$ (respectively $u_i$) in $S(\bar{V})$ (respectively $\bar{V}$), $i=2,...,n$.
If $\{\bar{u}_2,...,\bar{u}_n\}$ is linearly dependent over $F$, it violates the polynomial ring property of $F[\bar{z_2},...,\bar{z_n}]$.
It is now also clear that $\{v,u_2,...,u_n\}$ is a basis of $V$ over $F$.$\Box$

It is a well-known classical result that for a a finite $p$-group $G \subset GL(V)$, with finite $\dim_FV$ and $\Char F=p>0$, that there exists a basis of $V$,
putting $G$ simultaneously into a triangular form.
The next result is an extension of this in the non-linear setting.
The generators of the given polynomial ring are changed in a way reflecting the triangulating action of $G$ on $V$.

\begin{Proposition} \label {P2-9}
  Let $G \subset GL(V)$ be a finite $p$-group, where $n:=\dim _F V$ is finite, $p = \Char F > 0$ and $F=F^p$.
  Let $A =F[v^{p^{e_1}},x_2,...,x_n]$ be a graded polynomial subring of $S(V)$ (in $n$-variables) with $G \subset Aut _{\grade} (A)$. Assume:
  \begin{enumerate}
    \item $v \in A^G \cap V$;
    \item $x_i$ is homogenous, $\deg x_i = p^{e_i}$, $x_i - v_i^{p^{e_i}} \in v^{p^{e_1}}S(V)$, where $v_i \in V$, $i=2,...,n$, and $\{v,v_2,...,v_n\}$ is a basis of $V$;
    \item $e_1 \le e_2 \le \cdots \le e_n$.
  \end{enumerate}
  Then there exist homogenous elements $y_2,...,y_n \in A$, $\deg y_i=p^{e_i}$, $y_i-u_i^{p^{e_i}} \in v^{p^{e_1}}S(V)$, where $u_i \in V$, $i=2,...,n$, having the following properties:
  \begin{enumerate}
  \setcounter{enumi}{3}
    \item $(g-1)(y_2) \in F[v^{p^{e_1}}]_{+}$, $(g-1)(y_i)\in F[v^{p^{e_1}},y_2,...,y_{i-1}]_{+}$, $i>2$, for each $g\in G$;
    \item $A=F[v^{p^{e_1}},y_2,...,y_n]$;
    \item $\{v,u_2,...,u_n\}$ is a triangulating basis for the action of $G$ on $V$;
  \end{enumerate}
\end{Proposition}

\textbf{Proof:}
$\{x_2,...,x_n\}=\cup_{r=1}^l\{x_{i_{r-1}+1},...,x_{i_{r}}\}$, is a union of $l$ disjoint blocks, where $\deg x_i = \deg x_{i_r}=p^{e_{i_r}}$, for $i_{r-1}+1 \le i \le i_r$.
Set $x_{i_0}=v^{p^{e_1}}$.
Suppose $(4)$ holds for $i \le i_{r-1}:=s$ as well as $F[v^{p^{e_1}},x_2,...,x_s]=F[v^{p^{e_1}},y_2,...,y_s]$.
We shall now establish (4) for $i \le t:=i_r$, as well as $F[v^{p^{e_1}},x_2,...,x_t]=F[v^{p^{e_1}},y_2,...,y_t]$.
We therefore consider the variables $\{x_{s+1},...,x_t\}$, $\deg x_i=p^e$, $e=e_i$, $i_{r-1}+1\le i \le i_r$.
Since $g\in G$ preserves degrees we have:
\begin{equation}\label{eq2:4}
\begin{split}
  (g-1)(x_i) & =\sum_{j=s+1}^{t}b_{ij}(g)x_j+x,     \\
             & \text{ where } b_{ij}(g) \in F, \text{ } x \in F[v^{p^{e_1}},x_2,...,x_s]_{+},\text{ for each }g \in G.
\end{split}
\end{equation}
Similarly $(v^{p^{e_1}},x_2,...,x_s)$, the prime graded ideal generated by $\{v^{p^{e_1}},x_2,...,x_s\}$ in $F[v^{p^{e_1}},x_2,...,x_t]$, is $G$-stable.
Consequently $G$ acts by graded automorphisms on
$F[v^{p^{e_1}},x_2,...,x_t]/(v^{p^{e_l}},x_2,...,x_s) = F[\bar{x}_{s+1},...,\bar{x}_t]$, inducing by equality (\ref{eq2:4}) a linear action on $F\bar{x}_{s+1}+\cdots + F\bar{x}_t:=W$, a $(t-s)$-dimensional $F$-subspace of the polynomial ring $F[\bar{x}_{s+1},...,\bar{x}_t]$.
It implies, since $G$ is a $p$-group, the existence of a $W$-basis $\{z_{s+1},...,z_t\}$ which triangulates the action of $G$ on $W$. Thus:
\begin{equation}\label{eq2:5}
  (g-1)(z_{s+1}) = 0,\text{ } (g-1)(z_i)= \sum_{j=s+1}^{i-1}a_{ij}(g)z_j,\text{ } a_{ij}(g)\in F,\text{ }s+2\le i\le t.
\end{equation}
Since $\{\bar{x}_{s+1},...,\bar{x}_t\}$ is a basis of $W$ it follows that $z_i=\sum_{j=s+1}^{t}\lambda_{ij}\bar{x}_j$, $\lambda_{ij} \in F$, $i=s+1,...,t$ and $\lambda:=(\lambda_{ij})$ is an invertible matrix.
Set $y_i:=\sum_{j=s+1}^{t} \lambda_{ij}x_j$, $i=s+1,...,t$. Let $\lambda ^{-1}=(\mu _{ij})$. Then $x_i= \sum_{j=s+1}^{t}\mu _{ij}y_j$, $i=s+1,...,t$.
Therefore $F[v^{p^{e_1}},x_2,...,x_t]=F[v^{p^{e_1}},x_2,...,x_s][x_{s+1},...,x_t]=F[v^{p^{e_1}},y_2,...,y_s][y_{s+1},....,y_t]=F[v^{p^{e_1}},y_2,...,y_t]$.
Set $u_i:=\sum_{j=s+1}^{t}\lambda_{ij}^{\frac{1}{p^e}}v_j$, $i=s+1,...,t$.
Then $y_i-u_i^{p^e}=\sum_{j=s+1}^{t}\lambda_{ij}x_j -\sum_{j=s+1}^{t}\lambda_{ij}v_j^{p^e}=\sum_{j=s+1}^{t}\lambda_{ij}(x_j-v_j^{p^e})\in v^{p^{e_1}}S(V)$, $i=2,...,t$.

We now verify $(4)$,$(5)$ for $i\le t$.
Since $\bar{y}_i=z_i$, for $i=s+1,...,t$, it follows from the equalities (\ref{eq2:5}) that each member in $B:=\{(g-1)(y_{s+1}),(g-1)(y_i) -\sum_{j=s+1}^{i-1}a_{ij}(g)y_j$, $s+2 \le i \le t$, $g\in G\}$ is homogenous of degree $p^e$ (unless it is $0$), and is in $(v^{p^{e_1}},x_2,...,x_s)$.
This shows that each member in $B$, considered as sum of monomials in $F[v^{p^{e_1}},x_2,...,x_t]$, will not have $x_i$, with $i \ge s+1$ in any of its monomials because $\deg x_i=p^e$, for $i \ge s+1$ and each such monomial will contain at least one of $\{v^{p^{e_1}},x_2,...,x_s\}$.
Therefore each member of $B$ is in $F[v^{p^{e_1}},x_2,...,x_s]_{+}=F[v^{p^{e_1}},y_2,...,y_s]_{+}$.
This shows that $(4)$ holds for $i_{r-1}+1 \le i \le i_r$, as well.
Iterating this by induction will establish $(4)$, $(5)$.

We aim now at verifying $(6)$.
We shall firstly show that $\{v,u_2,...,u_n\}$ is a basis of $V$.
$v^{p^{e_1}}$ is a prime element in $A$, and hence $\mathfrak{q} = v^{p^{e_1}}A$ is a height one prime ideal in $A$.
Since $v^{p^{e_1}}S(V) \subseteq vS(V)$ it follows that $vS(V)\cap A = \mathfrak{p}$ is a prime ideal in $A$ with $\mathfrak{p} \supseteq \mathfrak{q}$.
Since $A/\mathfrak{p} =F[\bar{v}_2^{p^{e_2}},...,\bar{v}_n^{p^{e_n}}]$ (regarded as a subring of $S(V)/vS(V)$) it follows since $\{\bar{v}_2,...,\bar{v}_n\}$, is a basis of $V/Fv$, that $K.\dim A/\mathfrak{p} =n-1$.
Therefore $K.\dim A=n$ implies that $\height \mathfrak{p} = 1$ and $\mathfrak{q}=\mathfrak{p}$.
Hence $A/\mathfrak{p}=A/\mathfrak{q}=F[\bar{u}_2^{p^{e_2}},...,\bar{u}_n^{p^{e_n}}]$ is a polynomial ring in $(n-1)$ variables, implying that $\{\bar{u}_2,...,\bar{u}_n\}$ are linearly independent over $F$ and therefore a basis of $V/Fv$. Consequently $\{v,u_2,...,u_n\}$ is a basis of $V$.

We finally establish the triangulating property of $\{v,u_2,...,u_n\}$.
Recall that $y_i-u_i^{p^{e_i}} \in v^{p^{e_1}}S(V)$, $i=2,...,n$.
By $(4)$ $(g-1)(y_2)\in F[v]_{+}$, for each $g \in G$.
Hence $(g-1)(y_2)=\alpha v^{p^{e_2}}$, $\alpha \in F$.
Since $y_2-u_2^{p^{e_2}} \in v^{p^{e_1}}S(V)$ it follows that $(g-1)(u_2)^{p^{e_2}} \in vS(V)$, a prime ideal in $S(V)$.
Hence $(g-1)(u_2) \in vS(V)\cap V= Fv$, for each $g \in G$.

For $i>2$ we have by $(4)$ that $(g-1)(y_i)\in F[v,y_2,...,y_{i-1}]_{+}$
Hence $(g-1)(u_i)^{p^{e_i}}\in (v,u_2,...,u_{i-1})$ (meaning now the ideal in $S(V)$ generated by $\{v,u_2,...,u_{i-1}\}$).
Since the latter is a prime ideal then $(g-1)(u_i)\in (v,u_2,...,u_{i-1})\cap V=Fv+Fu_2+\cdots +Fu_{i-1}$, for each $g \in F$.$\Box$

We shall need the following old result of Nagata.
\begin{Proposition} \label{T2-10}
\cite[p. 143]{Na}. Let $R$ be a commutative Noetherian regular local ring and $\mathfrak{m}$ its unique maximal ideal.
Let $\mathfrak{p}$ be a prime ideal in $R$.
Then $\mathfrak{p}^{(e)} \subseteq \mathfrak{m}^e$, for every $e \in \mathbb{N}$, where $\mathfrak{p}^{(e)}:=\mathfrak{p}_{\mathfrak{p}}^e \cap R$ is the $e$-th symbolic power of $\mathfrak{p}$.
\end{Proposition}

The next result provides an hypersurface section reduction which will be used repeatedly in Theorem \ref{T2-12}.

\begin{Proposition} \label {P2-11}
  Let $A=\oplus_{n \ge 0} A_n$ be an $\mathbb{N}$-graded polynomial ring over a field $F =A_0$ and $G\subset Aut_{\grade}(A)$ a finite subgroup of graded automorphisms.
  Let $x \in A^G$ be a homogenous element of $A$.
  Set $H:=\ker(G \rightarrow G|_{A/ xA})$, $\mathfrak{m}:= A_{+}$, the irrelevant maximal ideal of $A$ and $\mathfrak{n}:=A_{+}^G$.
  Suppose:
  \begin{enumerate}
    \item $x \in \mathfrak{m} - \mathfrak{m}^2$;
    \item $G/H$ acts faithfully on   $A^H/xA^H$;
    \item $A_{\mathfrak{n}}^G$ satisfies Serre's $S_{m+1}$ and $R_m$ conditions with $m \ge 2$.
  \end{enumerate}
  Then $A^G/(x)=(A^H/xA^H)^{G/H}$, and $[A^G/(x)]_{\mathfrak{n}/(x)}$ satisfies $S_m$ and $R_{m-1}$.
\end{Proposition}

\textbf{Proof:}
By $(1)$ and \cite[Theorem 161]{Ka} $A_\mathfrak{m}/xA_\mathfrak{m}=[A/xA]_{\bar{\mathfrak{m}}}$ is a regular local ring, where $\bar{\mathfrak{m}}=\mathfrak{m}/(x)$.
Hence it is a domain implying the same for $A/xA$, so $x$ is a prime element in $A$.
Since $xA\cap A^H=xA^H$ and $xA\cap A^G = xA^G := (x)$, it follows that $x$ is also a prime element in $A^H$ and $A^G$.

We next show that $[A^G/(x)]_{\bar{\mathfrak{n}}}$ satisfies $R_{m-1}$ where $\bar{\mathfrak{n}}=\mathfrak{n}/(x)$.
Let $\mathfrak{p} \subset \mathfrak{n}$ be a prime ideal in $A^G$ with $x \in \mathfrak{p}$, $\height(\mathfrak{p})=m$.
Hence $A^G_{\mathfrak{p}}$ is regular.
Let $\{ \mathfrak{p}_1,...,\mathfrak{p}_k\}$ be the set of prime ideals in $A$ lying over $\mathfrak{p}$.
Since $\mathfrak{m}$ is the unique maximal ideal in $A$ lying over $\mathfrak{n}$, it follows by "going up" that $\mathfrak{p}_i \subset \mathfrak{m}$, $i=1,...,k$.
We have:\\
$\mathfrak{p}_i^{(2)}=
\mathfrak{p}_i^2 A_{\mathfrak{p}_i} \cap A
= \mathfrak{p}_{i \mathfrak{m}}^2 (A_{\mathfrak{m}})_{\mathfrak{p}_{im}}\cap A_{\mathfrak{m}} \cap A
= \mathfrak{p}_{i \mathfrak{m}}^{(2)} \cap A \subset \mathfrak{m}_{\mathfrak{m}}^2 \cap A
= \mathfrak{m} ^2$,
where the inclusion is due to Proposition \ref{T2-10} and the last equality is by the maximality of $\mathfrak{m}$.
Consequently by $(1)$ $x \notin \mathfrak{p}_i^{(2)}$.
Since
$\mathfrak{p}^{(2)}
=\mathfrak{p}^2A_{\mathfrak{p}}^G \cap A^G \subset \mathfrak{p}_i^2 A_{\mathfrak{p}_i} \cap A
=\mathfrak{p}_i^{(2)}$,
it follows that $x \notin \mathfrak{p}^{(2)}$.
Hence $x \in \mathfrak{p}_{\mathfrak{p}} - \mathfrak{p}_\mathfrak{p}^2$ and therefore, by the regularity at $A_{\mathfrak{p}}^G$ and by \cite[Thm. 161]{Ka}, we have that $A_{\mathfrak{p}}^G/xA_{\mathfrak{p}}^G =[A^G/(x)]_{\bar{\mathfrak{p}}}$ is regular where $\bar{\mathfrak{p}}=\mathfrak{p}/(x)$.
Moreover $K.\dim [A^G/(x)]_{\bar{\mathfrak{p}}}=m-1$.
So $[A^G/(x)]_{\bar{\mathfrak{n}}}$ satisfies $R_{m-1}$.

Similarly suppose $\bar{\mathfrak{p}} =\mathfrak{p}/(x)$ is a prime ideal in $A^G/(x)$ with $\bar{\mathfrak{p}} \subset \bar{\mathfrak{n}}$.
Since $A_{\mathfrak{n}}^G$ satisfies $S_{m+1}$ it follows that $\depth \mathfrak{p}_{\mathfrak{p}}=\depth (\mathfrak{p}_{\mathfrak{n}})_{\mathfrak{p}_{\mathfrak{n}}} \ge min\{m+1, \height \mathfrak{p}_{\mathfrak{n}}\} = min\{m+1, \height \mathfrak{p}\}$.
Since $(x) \subset \mathfrak{p}$, $x$ can be completed to a regular sequence inside $\mathfrak{p}_{\mathfrak{p}}$, implying that $\depth \bar{\mathfrak{p}}_{\bar{\mathfrak{p}}} = \depth \mathfrak{p}_{\mathfrak{p}}-1$.
Since $\height \bar{\mathfrak{p}} = \height \mathfrak{p} - 1$, the inequality $\depth \bar{\mathfrak{p}}_{\bar{\mathfrak{p}}} \ge min\{m, \height \bar{\mathfrak{p}} \}$ clearly holds and the $S_m$ property of $(A^G/(x))_{\bar{\mathfrak{n}}}$ is established.

Consequently, since $m \ge 2$, $[A^G/(x)]_{\bar{\mathfrak{n}}}$ satisfies $S_2$ and $R_1$, implying that it is normal.

Now by $(2)$ and Galois theory we have that $[Q(A^H/xA^H):Q(A^H/xA^H)^{G/H}]=|G/H|$.
Since $G$ is finite, $Q(A^H/xA^H)^{G/H}=Q((A^H/xA^H)^{G/H})$ and we therefore have:
\begin{equation}\label{eq2:7}
  [Q(A^H/xA^H):Q((A^H/xA^H)^{G/H}]=|G/H|.
\end{equation}
Similarly $(A^H)^{G/H}=A^G$, $(Q(A^H))^{G/H}=Q(A^G)$ and $[Q(A)^H:Q(A^G)]=|G/H|$.
Since $A_{(x)}^G$ is a DVR it follows that $A_{(x)}^H$ is free of rank $|G/H|$ over $A_{(x)}^G$.
Consequently:
\begin{equation*} \label{e2:8.1}
 [A^H_{(x)}/xA^H_{(x)}:A^G_{(x)}/(x)_{(x)}]=|G/H|, \text{implying: }
\end{equation*}
\begin{equation}\label{eq2:8}
  [Q(A^H/xA^H):Q(A^G/(x))]=|G/H|.
\end{equation}
Since $A^G/(x) \subseteq (A^H/xA^H)^{G/H}$, it follows from equalities (\ref{eq2:7}) and (\ref{eq2:8}) that $Q(A^G/(x)) =Q((A^H/xA^H)^{G/H})$.
Now since $[A^G/(x)]_{\bar{\mathfrak{n}}}$ is normal and $(A^H/xA^H)^{G/H}$ is finite over $A^G/(x)$, it follows that:
\begin{equation}\label{eq2:9}
  [A^G/(x)]_{\bar{\mathfrak{n}}}=[(A^H/xA^H)^{G/H}]_{\bar{\mathfrak{n}}}.
\end{equation}
Finally since $A^G/(x) \subseteq (A^H/xA^H)^{G/H}$ is a finite extension and both are $\mathbb{N}$-graded with $F$ as the $0$-component, it follows from equality (\ref{eq2:9}) that there exists $y= \alpha +z$, $\alpha \in F - \{0\}$, $z \in \bar{\mathfrak{n
}}$, with $y[A^H/xA^H] \subseteq A^G/(x)$.
Hence $[A^H/xA^H]=\bar{\mathfrak{n}}[A^H/xA^H]+[A^G/(x)]$.
So by the graded version of Nakayama's lemma (e.g. \cite[1.5.24]{BH}) $(A^H/xA^H)^{G/H}=A^G/(x)$.$\Box$

The following is a useful criterion on a given system of $n$($=\dim _F V$) homogenous elements of $S(V)^G$.

\begin{Proposition} \label{P2-11.5}
\cite[Thm. 3.9.4]{DK}. Let $G \subset GL(V)$ is a finite group, $\dim _FV=n$ and $S(V)=F[x_1,...,x_n]$.
Let $f_1,...,f_n \in S(V)^G$ be homogenous elements.
Then the following statements are equivalent:
\begin{enumerate}
\item $S(V)^G=F[f_1,...,f_n]$ (in particular $S(V)^G$ is a polynomial ring);
\item $\{f_1,....,f_n\}$ are algebraically independent over $F$ and $\prod_{i=1}^{n} \deg (f_i)=|G|$;
\item The Jacobian determinant $\det (\frac{\partial f_i}{\partial x_j})$ is non-zero and $\prod_{i=1}^{n} \deg (f_i)=|G|$;
\item $\{f_1,...,f_n\}$ is a system of parameters (s.o.p.) and $\prod_{i=1}^{n} \deg (f_i)=|G|$.
\end{enumerate}
\end{Proposition}
\textbf{Proof:} The equivalence of items $(1)$, $(2)$ and $(3)$ are due to \cite{K-Cal}. The equivalence of items (1) and (4) is standard and follows from (e.g. \cite[Prop. 3.11.2]{DK}).

The following is the proof of Theorem B in case $F$ is perfect.

\begin{Theorem} \label {T2-12}
  Let $G$ be a finite $p$-group, where $n:=\dim _F V \ge 4$, $\Char F = p$, and $F=F^p$. Then $S(V)^G$ is a polynomial ring if and only if:
    \begin{enumerate}
    \item $S(V)^{G_U}$ is a polynomial ring, for each subspace $U \subset V^*$ with $\dim _F U=2$, and
    \item $S(V)^G$ is Cohen-Macaulay.
    \end{enumerate}
  \end{Theorem}

\textbf{Proof:}
By condition $(1)$ and Lemma \ref{L2-4} $S(V)^G$ satisfies Serre's $R_{n-2}$ condition.
Hence by Lemma \ref{L2-5} and assumption $(2)$ $G$ is generated by transvections.
$G$ is $p$-group so by a classical result \cite[Theorem C, p.100]{Kap} there is a $G$-triangulating basis $\{v_1,y_2^{(1)},...,y_n^{(1)}\}$ of $V$, for the action of $G$ on $V$, with $g(v_1)=v_1$,
and $(g-1)(y_j^{(1)}) \in Fv_1+Fy_2^{(1)}+\cdots +Fy_{j-1}^{(1)}$, $j>2$, $(g-1)(y_2^{(1)}) \in Fv_1$,
$\forall g \in G$.

Set $A_1=S(V))$, $x_1=v_1$, $G_1=G$, $V_1=V$.
For $i=1$ the next (I),(II),(III),(IV),(V) hold by the previous paragraph.
Given by induction the existence of a quadruple $(A_i, x_i,V_i, G_i)$, $i \ge 2$, with the following properties:
   \begin{enumerate}[label=(\Roman*)]
   \item $V_i:=V/Fv_1+\cdots +Fv_{i-1}$, $\dim _F V_i=n-i+1$, $\{v_1,...,v_{i-1} \} \subseteq V$;
   \item $G_i:=G/G_{(Fv_1+\cdots +Fv_{i-1})^{\perp}}$, $G_i$ acts faithfully on $V_i$, where this action is induced from the action of $G$ on $V$;
   \item $A_i:=F[x_i,y_{i+1}^{(i)},...,y_n^{(i)}] \subset S(V_i)$ is a polynomial ring in $(n-i+1)$ homogenous elements, $x_i=\bar{v}_i^{p^{e_i}}$, $y_j^{(i)}=\bar{u}_j^{p^{f_j}}$, $j=i+1,...,n$, where $\{\bar{v}_i,\bar{u}_{i+1},...,\bar{u}_n\}$ is a triangulating basis for the action of $G_i$ on $V_i$ and $e_i\le f_{i+1}\le \cdots \le f_n$;
   \item $(g-1)(x_i)=0$, $(g-1)(y_{i+1}^{(i)})\in F[x_i]_{+}$, $(g-1)(y_j^{(i)}) \in F[x_i,y_{i+1}^{(i)},...,y_{j-1}^{(i)}]_{+}$ for each $g \in G_i$ and $j=i+2,...,n$;
   \item $A_{i-1}^{G_{i-1}}/(x_{i-1})=A_i^{G_i}$, $[A_i^{G_i}]_{\mathfrak{n}}$ is Cohen-Macaulay and satisfies Serre's $R_{n-i-1}$ condition, where $\mathfrak{n}:= (A_i^{G_i})_{+}$.
   \end{enumerate}

We shall now proceed to construct the $(i+1)$-th quadruple with the above properties.

\textbf{Step 1:} \underline{$G$ acts in a triangular way on $Fv_1+\cdots + Fv_i$, and consequently $(v_1,...,v_i)$ is a}
\\\underline{$G$-stable prime ideal in $S(V)$.}

We have, from properties (III) and (IV), which are valid by induction for $j \le i$, that $(g-1)(\bar{v}_j)^{p^{e_j}}=(g-1)(x_j)=0$, $\forall g \in G_j$.
Therefore, using $S(V_j)=S(V)/(v_1,...,v_{j-1})$, $(g-1)(v_j)^{p^{e_j}} \in (v_1,...,v_{j-1})$, $\forall g \in G$.
So since the latter is a prime ideal in $S(V)$, we get by degree considerations that $(g-1)(v_j)\in (v_1,...,v_{j-1})\cap V=Fv_1+\cdots +Fv_{j-1}$, $\forall g \in G$, and $2 \le j \le i$.
Consequently $(v_1,...,v_i)$ is a $G$-stable ideal in $S(V)$.

\textbf{Step 2:} \underline{$x_iA_i=\bar{v}_iS(V_i)\cap A_i$, and consequently $A_i/x_iA_i \subset S(V_{i+1})$, where $V_{i+1}:=$ }
\\\underline{$V/Fv_1+\cdots +Fv_{i-1}+Fv_i$.}

By the polynomial property in (III) we have that $\mathfrak{q}:= x_iA_i$ is a prime ideal in $A_i$ and $\height (\mathfrak{q})=1$.
Since $\bar{v}_iS(V_i)$ is a (height one) prime ideal in $S(V_i)$, it follows that $\mathfrak{p}:=\bar{v}_iS(V_i)\cap A_i$ is a prime ideal in $A_i$ with $\mathfrak{p}\supseteq \mathfrak{q}$.
Since $\{\bar{v}_i,\bar{u}_{i+1},...,\bar{u}_n\}$ is a basis of $V_i=V/Fv_1+\cdots +Fv_{i-1}$, it follows that the images $\{\tilde{u}_{i+1},...,\tilde{u}_n\}$ form a basis of $V_{i+1}$, so they are generators of $S(V_{i+1})=S(V_i)/\bar{v}_iS(V_i)=S(V)/(v_1,...,v_i)$.
Therefore $A_i/\mathfrak{p}=F[\bar{y}_{i+1}^{(i)},...,\bar{y}_n^{(i)}]=F[\tilde{u}_{i+1}^{p^{f_{i+1}}},...,\tilde{u}_n^{p^{f_n}}]$ is a polynomial ring in $(n-i)$ variables, implying since $K.\dim A_i=(n-i+1)$, that $\height (\mathfrak{p})=1$ and hence $\mathfrak{p}=\mathfrak{q}$.

\textbf{Step 3:} \underline{$H_i=G_{(Fv_1+\cdots +Fv_{i})^{\perp}} /G_{(Fv_1+\cdots +Fv_{i-1})^{\perp}}$, where $H_i:=\{g \in G_i|(g-1)(A_i)\subseteq x_iA_i\}$.}

Let $g \in H_i$, then $(g-1)(y_j^{(i)})\in x_iA_i$, $j=i+1,...,n$.
Hence $(g-1)(y_j^{(i)})\in x_iS(V_i) \subseteq \bar{v}_iS(V_i)=(v_1,...,v_i)/(v_1,...,v_{i-1})$.
Write $g=g'G_{(Fv_1+\cdots Fv_{i-1})^{\perp}}$, $g' \in G$.
It follows that $[(g'-1)(u_j)]^{p^{f_j}} \in (v_1,...,v_i)$.
Since the latter ia a prime ideal in $S(V)$, we get that $(g'-1)(u_j)\in (v_1,...,v_i)\cap V=Fv_1+\cdots +Fv_i$, $j=i+1,...,n$.
Also $0=(g-1)(x_i)=[(g-1)(\bar{v}_i)]^{p^{e_i}}$, implying that $(g-1)(\bar{v}_i)=0$, hence $(g'-1)(v_i)\in Fv_1+\cdots +Fv_{i-1}$.
Consequently $(g'-1)(V)\subseteq Fv_1+\cdots +Fv_i$.
Equivalently $g' \in G_{(Fv_1+\cdots +Fv_i)^{\perp}}$ and hence $g \in G_{(Fv_1+\cdots +Fv_i)^{\perp}}/G_{(Fv_1+\cdots +Fv_{i-1})^{\perp}}$, and $H_i \subseteq G_{(Fv_1+\cdots +Fv_i)^{\perp}}/G_{(Fv_1+\cdots +Fv_{i-1})^{\perp}}$. The reversed inclusion is clear.

\textbf{Step 4:} \underline{$A_i^{H_i}$ is a polynomial ring.}

Recall that $A_i=F[x_i,y_{i+1}^{(i)},...,y_n^{(i)}] \subset S(V_i)$ is a polynomial ring in $(n-i+1)$ homogenous elements, and $x_i=\bar{v}_i^{p^{e_i}}$, $y_j^{(i)}=\bar{u}_j^{p^{f_j}}$, $j=i+1,...,n$, $e_i\le f_{i+1} \le \cdots \le f_n$.
We have from the previous paragraph that $(g'-1)(u_j)\in Fv_1+...+Fv_i$.
Consequently $(g-1)(\bar{u}_j)\in F\bar{v}_i$ and $g(\bar{v}_i)=\bar{v}_i$, $\forall g\in H_i$, $j=i+1,...,n$.
It follows from Corollary \ref{C2-8} that $A_i^{H_i}=F[x_i,z_{i+1},...,z_n]$ is a polynomial ring in $(n-i+1)$ homogenous elements, $\deg z_j=p^{h_j}$, $z_j-\bar{w}_j^{p^{h_j}} \in \bar{v}_i^{p^{e_i}}S(V_i)$, for $j=i+1,...,n$, where $e_i \le h_{i+1}\le \cdots \le h_n$ and $\{\bar{v}_i,\bar{w}_{i+1},...,\bar{w}_n\}$ is a basis of $V_i$.

\textbf{Step 5:} \underline{$G_i/H_i$ acts faithfully on  $A_i^{H_i}/x_iA_i^{H_i}$.}

Indeed if $g \in G_i$, with $(g-1)(A_i^{H_i}) \subseteq x_iA_i^{H_i}=\bar{v}_i^{p^{e_i}}A_i^{H_i} \subseteq \bar{v}_i^{p^{e_i}}S(V_i)$, then
$(g-1)(\bar{w}_j)^{p^{h_j}}=(g-1)(\bar{w}_j^{p^{h_j}}) \in  \bar{v}_iS(V_i)$.
Since the latter is a prime in $S(V_i)$, it follows that $(g-1)(\bar{w}_j)\in \bar{v}_iS(V_i)\cap V_i=F\bar{v}_i$, for $j=i+1,...,n$.
Since $V_i =\Span_F \{\bar{v}_i,\bar{w}_{i+1},...,\bar{w}_n\}$ it follows that $(g-1)(V_i)\subseteq F\bar{v}_i$.
Equivalently if $g=g'G_{(Fv_1+\cdots +Fv_{i-1})^{\perp}}$ then $(g'-1)(V) \subseteq Fv_1+\cdots +Fv_i$, and hence $g'\in G_{(Fv_1+\cdots +Fv_i)^{\perp}}$.
Consequently $g \in H_i$ as claimed.

\textbf{Step 6:} \underline{Changing the generators of $A_i^{H_i}$.}

We shall now apply Proposition \ref{P2-9} with $G_i$ acting on $A_i^{H_i}$. Then $A_i^{H_i}=F[x_i,m_{i+1},...,m_n]$, $m_j-\bar{t}_j^{p^{h_j}}\in \bar{v}_i ^{p^{e_i}}S(V_i)$ for $j=i+1,...,n$, $\{\bar{v}_i,\bar{t}_{i+1},...,\bar{t}_n\}$ is a triangulating basis for the action of $G_i$ on $V_i$, $\deg(m_j)=p^{h_j}$, and $m_j$ is a homogenous element in $S(V_i)$, $j=i+1,...,n$.
Moreover $(g-1)(m_{i+1})\in F[\bar{v}_i^{p^{e_i}}]_{+}=F[x_i]_{+}$, $(g-1)(m_j)\in F[x_i,m_{i+1},...,m_{j-1}]_{+}$, for $j>i+1$, $g\in G_i$.

\textbf{Step 7:} \underline{$(A_{i+1},x_{i+1},V_{i+1},G_{i+1})$ and its properties.}

Set $A_{i+1}:=A_i^{H_i}/x_iA_i^{H_i}$, $V_{i+1}:=V/Fv_1+\cdots +Fv_i$ (as in Step 2), $G_{i+1}:=G_i/H_i=G/G_{(Fv_1+\cdots Fv_i)^{\perp}}$ and $x_{i+1}:=\bar{m}_{i+1}$, the image of $m_{i+1}$ in $A_{i+1} \subset S(V_{i+1})$, where the last inclusion is by Step 2.

We next verify property $(II)$ for $(A_{i+1},x_{i+1},V_{i+1},G_{i+1})$.
By Step 1, the action of $G$ on $V_{i+1}$ is well defined and linear, so inducing a linear action of $G_{i+1}=G/G_{(Fv_1+\cdots Fv_i)^{\perp}}$ on $V_{i+1}$. We now show that this action is faithful.
Say $(g-1)(V_{i+1})=0$, for some $g \in G_{i+1}=G/H_i$.
Then $(g-1)(S(V_{i+1}))=0$.
Hence $(g-1)(A_i^{H_i}/x_iA_i^{H_i})=0$ and $g=1$ follows from Step 5.

Set $e_{i+1}:=h_{i+1}$, $v_{i+1}:=t_{i+1}$, $y_j^{(i+1)}:=\bar{m}_j$, the image of $m_j$ in $A_{i+1}$, $j=i+2,...,n$.
Since $A_i^{H_i}=F[x_i,m_{i+1},...,m_n]$ is a polynomial ring, it follows that $A_{i+1}=A_i^{H_i}/x_iA_i^{H_i}$ is a polynomial ring in $(n-i)$ variables.
This implies that $A_{i+1}=F[x_{i+1}, y_{i+2}^{(i+1)},...,y_n^{i+1}]$, $x_{i+1}=\bar{v}_{i+1}^{p^{e_{i+1}}}$, $y_j^{(i+1)}=\bar{t}_j^{p^{h_j}}$, $j=i+2,...,n$ (by abuse of notation we still denote by $\bar{t}_j$ the image of $t_j$  in $V_{i+1})$.
It now follows from Step 6 that $(g-1)(x_{i+1})=0$, $(g-1)(y_{i+2}^{(i+1)})=(g-1)(\bar{m}_{i+2})\in F[\bar{m}_{i+1}]_{+}=F[x_{i+1}]_{+}$, $(g-1)(y_j^{(i+1)})\in F[x_{i+1},y_{i+2}^{(i+1)},...,j_{j-1}^{(i+1)}]_{+}$, $j>i+2$, $\forall g \in G_{i+1}$.
It also follows, as in the proof of Proposition \ref{P2-9}(6),  that $\{\bar{v}_{i+1},\bar{t}_{i+2},...,\bar{t}_n\}$ is a triangulating basis for the action of $G_{i+1}$ in $V_{i+1}$.

This settles properties $(I)$, $(II)$, $(III)$ and $(IV)$ of $(A_{i+1},x_{i+1},V_{i+1},G_{i+1})$.

Now we verify property $(V)$ for the quadruple $(A_{i+1},x_{i+1},V_{i+1},G_{i+1})$.

Set $\mathfrak{n}:=[A_i^{G_i}]_{+}=$ the irrelevant maximal ideal of $A_i^{G_i}$, $\mathfrak{m}:=[A_i]_{+}=$ the irrelevant maximal ideal of $A_i$.
Since $x_i\in \mathfrak{m}-\mathfrak{m}^2$ (by the fact that it is a minimal generator  of $A_i$), it is also a prime element of the polynomial ring $A_i$ as well as being a $G_i$-invariant.
By Step 5 we have that $G_i/H_i$ acts faithfully on $A_i^{H_i}/x_iA_i^{H_i}$.
Also by property ($V)$ for $(A_{i},x_{i},V_{i},G_{i})$, $[A_i^{G_i}]_{\mathfrak{n}}$ satisfies Serre's $R_{n-i-1}$ condition , as well as being Cohen-Macaulay. This implies by Proposition \ref{P2-11} that:\\
$A_i^{G_i}/(x_i)=(A_i^{H_i}/x_iA_i^{H_i})^{G_i/H_i} = A_{i+1}^{G_{i+1}}$, and $[A_i^{G_i}/(x_i)]_{\mathfrak{n}/(x_i)}=[A_{i+1}^{G_{i+1}}]_{\mathfrak{n}/(x_i)}$ satisfies Serre's $R_{n-i-2}$ condition.
It is also a Cohen-Macaulay ring, so property $(V)$ is verified for the quadruple $(A_{i+1},x_{i+1}, V_{i+1},G_{i+1})$.

This settles the inductive step.

\textbf{Step 8:} \underline{The chain of hypersurface sections of $S(V)^G$.}

Due to the requirement of $m \ge 2$ in Serre's condition $R_m$, in Proposition \ref{P2-11}, this procedure can be carried only for $i \le n-3$.

Consequently we have the following chain of equalities:
\begin{equation} \label{eq2:16}
S(V)^G/(v)=A_1^{G_1}/(x_1)=A_2^{G_2}, \text{ } A_2^{G_2}/(x_2)=A_3^{G_3},..., A_{n-3}^{G{n-3}}/(x_{n-3})=A_{n-2}^{G_{n-2}}.
\end{equation}
Observe that $A_{n-2}=F[\bar{w}_{n-2}^{p^f}, \bar{w}_{n-1}^{p^g},\bar{w}_n^{p^h}] \subset S(V_{n-2})$, $V_{n-2}=V/Fv_1+\cdots +v_{n-3}$, $G_{n-2}=G/G_{(Fv_1+\cdots Fv_{n-3})^{\perp}}$, $f \le g \le h$ and $\{\bar{w}_{n-2}$, $\bar{w}_{n-1},\bar{w}_n\}$ is a triangulating basis for the faithful, linear action of $G_{n-2}$ on $V_{n-2}= F\bar{w}_{n-2}+F\bar{w}_{n-1}+F\bar{w}_n$.

\textbf{Step 9:} \underline{The chain of hypersurface sections of $S(V)^T$, where $T:=G_{(Fv_1+\cdots +Fv_{n-3})^{\perp}}$.}

Let $U:=(Fv_1+...+Fv_{n-3})^{\perp}$. $\dim_FU=n-(n-3)=3$ and therefore by assumption $(1)$ and Lemma \ref{L2-4}, $S(V)^{G_U}$ is a polynomial ring.
For each $1\le i\le n-3$, we have $G_{(Fv_1+\cdots +Fv_i)^{\perp}}=\{g \in G|(g-1)(V)\subseteq Fv_1+\cdots +Fv_i\} \subseteq \{g\in G|(g-1)(V)\subseteq Fv_1+\cdots +Fv_{n-3}\}=T$.
Consequently $G_{(Fv_1+\cdots +Fv_i)^{\perp}}=T_{(Fv_1+\cdots +Fv_i)^{\perp}}$, and therefore for $i  \ge 2$ $H_i=G_{(Fv_1+\cdots +Fv_i)^{\perp}}/G_{(Fv_1+\cdots +Fv_{i-1})^{\perp}}=T_{(Fv_1+\cdots +Fv_i)^{\perp}}/T_{(Fv_1+\cdots +Fv_{i-1})^{\perp}}$. Also for $i=1$, we have $H_1=G_{(Fv_1)^{\perp}}=T_{(Fv_1)^{\perp}}$.

Set $T_i:=T/T_{(Fv_1+\cdots +Fv_{i-1})^{\perp}}$, $i=2,...,n-2$.
Then since $S(V)^T$ is a polynomial ring we can similarly apply Proposition \ref{P2-11} with the same $A_i$, $H_i$, $x_i$, for each $i$, and get, as in Step 8, the following chain of equalities.
\begin{equation}\label{eq2:17}
  S(V)^T/(v_1)=S(V)^T/(x_1)=A_2^{T_2}, \text{ }A_2^{T_2}/(x_2) =A_3^{T_3},...,A_{n-3}^{T_{n-3}}=A_{n-2}^{T_{n-2}}.
\end{equation}
But $T_{n-2}=T/T_{(Fv_1+\cdots +Fv_{n-3})^{\perp}}=1$, implying that $A_{n-2}^{T_{n-2}}=A_{n-2}=F[\bar{w}_{n-2}^{p^f}, \bar{w}_{n-1}^{p^g},\bar{w}_n^{p^h}]  \subseteq S(V_{n-2})$. Also $f \le g \le h$ and $\{\bar{w}_{n-2},\bar{w}_{n-1},\bar{w}_n\}$ is a triangulating basis for the action of $G_{n-2}$ on $V_{n-2}=F\bar{w}_{n-2}+F\bar{w}_{n-1}+F\bar{w}_n$.

Let $\{x'_1,...,x'_{n-3},x'_{n-2},x'_{n-1},x'_n\}$ be homogenous preimages of $\{x_1,...,x_{n-3},\bar{w}_{n-2}^{p^f},\bar{w}_{n-1}^{p^g},\bar{w}_n^{p^h}\}$ in $S(V)^T$.
They exist by equalities (\ref{eq2:17}) and they generate $S(V)^T$.
Consequently by Proposition \ref{P2-11.5}:
\begin{equation}\label{eq2:14}
|T|=\prod_{i=1}^{n}(\deg x'_i)=\prod_{i=1}^{n-3}(\deg x_i)p^fp^gp^h.
\end{equation}

\textbf{Step 10:} \underline{The conclusion of the proof.}

We may clearly assume that $G \neq T$. $G_{n-2}=G/G_{(Fv_1+\cdots +Fv_{n-3})^{\perp}}=G/T$ acts linearly on $V_{n-2}=F\bar{w}_{n-2}+F\bar{w}_{n-1}+R\bar{w}_n$, and is generated by transvections (since $G$ is generated by transvections on $V$ by Lemma \ref{L2-5}, and by Lemma \ref{L2-6}, this property is inherited by $G_i$ on $V_i$, for each $i$).
Moreover $\{\bar{w}_{n-2},\bar{w}_{n-1},\bar{w}_n\}$ is a triangulating basis for the action of $G_{n-2}$ on $V_{n-2}$.
Therefore (as e.g. in Proposition \ref{P3-2}) $G_{n-2}$ is represented, with respect to this basis, by matrices in
$\begin{bsmallmatrix}
    1 & F & F \\
    0 & 1 & F \\
    0 & 0 & 1
  \end{bsmallmatrix}$.
We shall split the discussion into three cases.

\textbf{Case 1:} $G_{n-2} \subseteq \begin{bsmallmatrix}
    1 & 0 & F \\
    0 & 1 & F \\
    0 & 0 & 1
  \end{bsmallmatrix}$.\\
In this case we have that $G_{n-2}=H_{n-2}$, implying by Corollary \ref{C2-8}
that:
\begin{equation*}\label{eq2:18}
A_{n-2}^{G_{n-2}}=F[\bar{w}_{n-2}^{p^f},\bar{w}_{n-1}^{p^g},\bar{w}_n^{p^h}]^{H_{n-2}}=F[x_{n-2},x_{n-1},x_n].
\end{equation*}
Therefore if $\{x'_1,...,x'_{n-3},x'_{n-2},x'_{n-1},x'_n\}$ are homogenous preimages in $S(V)^G$ of
\\$\{x_1,...,x_{n-3},x_{n-2},x_{n-1},x_n\}$ which exist by the chain of equalities (\ref{eq2:16}), then they generate $S(V)^G$ and consequently $S(V)^G$ is a polynomial ring.

\textbf{Case 2:} $G_{n-2} \subset \begin{bsmallmatrix}
    1 & F & F \\
    0 & 1 & 0 \\
    0 & 0 & 1
  \end{bsmallmatrix}$.\\
Hence $g(\bar{w}_{n-2})=\bar{w}_{n-2}$, $g(\bar{w}_{n-1}=\bar{w}_{n-1}$ for each $g \in G_{n-2}$.
Consequently by Proposition \ref{P2-6.5} $F[\bar{w}_{n-2},\bar{w}_{n-1},\bar{w}_n]^{G_{n-2}}= F[\bar{w}_{n-2},\bar{w}_{n-1},a]$ is a polynomial ring where $\deg a=|G_{n-2}|=|G/T|$.
Therefore using $f \le g \le h$, $a^{p^h}\in F[\bar{w}_{n-2}^{p^f},\bar{w}_{n-1}^{p^g},\bar{w}_n^{p^h}]^{G_{n-2}}=A_{n-2}^{G_{n-2}}$.
Clearly $\{\bar{w}_{n-2}^{p^f},\bar{w}_{n-1}^{p^g},a^{p^h}\}$ is a system of parameters of $A_{n-2}^{G_{n-2}}$.
By the chain of equalities (\ref{eq2:16}) we can pick in $S(V)^G$ homogenous preimages $\{x'_1,...,x'_{n-3}, x'_{n-2},x'_{n-1},x'_n\}$ of
$\{x_1,...,x_{n-3}, \bar{w}_{n-2}^{p^f},\bar{w}_{n-1}^{p^g},a^{p^h}\}$.
These $n$ preimages form a system of homogenous parameters in $S(V)^G$, and by equality (\ref{eq2:14}):
\begin{equation*}\label{eq2:19}
\begin{split}
\prod_{i=1}^{n}(\deg x'_i)
& = \prod_{i=1}^{n-3}(\deg x'_i) p^f  p^g \deg a^{p^h}\\
& = \prod_{i=1}^{n-3}(\deg x_i) p^f p^g (\deg a)p^h
= |T| (\deg a)
= |T||G/T|
= |G|.
\end{split}
\end{equation*}
Consequently by Proposition \ref{P2-11.5} $S(V)^G=F[x'_1,...,x'_{n-3},x'_{n-2},x'_{n-1},x'_n]$ is a polynomial ring.

\textbf{Case 3:} $G_{n-2} \not \subset \begin{bsmallmatrix}
    1 & 0 & F \\
    0 & 1 & F \\
    0 & 0 & 1
  \end{bsmallmatrix}$,
$\begin{bsmallmatrix}
    1 & F & F \\
    0 & 1 & 0 \\
    0 & 0 & 1
  \end{bsmallmatrix}$.
\\
Therefore by Theorem \ref{T3-3}, $F[\bar{w}_{n-2},\bar{w}_{n-1},\bar{w}_n]^{G_{n-2}}=F[\bar{w}_{n-2},a,b]$, where $a\in F[\bar{w}_{n-2},\bar{w}_{n-1}]^{G_{n-2}}$, and hence $(\deg a)(\deg b)=|G_{n-2}|$.
Since $f \le g \le h$, then $\bar{w}_{n-2}^{p^f}$, $a^{p^g} \in F[\bar{w}_{n-2}^{p^f},\bar{w}_{n-1}^{p^g}]^{G_{n-2}}$,
$b^{p^h} \in F[\bar{w}_{n-2}^{p^f},\bar{w}_{n-1}^{p^g},\bar{w}_n^{p^h}]^{G_{n-2}}$.
In particular $\{\bar{w}_{n-2}^{p^f},a^{p^g},b^{p^h}\}$ is a system of parameters (s.o.p.) of $A_{n-2}^{G_{n-2}}=F[\bar{w}_{n-2}^{p^f},\bar{w}_{n-1}^{p^g},\bar{w}_n^{p^h}]^{G_{n-2}}$.
Therefore by the chain of equalities (\ref{eq2:16}) there are $n$ homogenous preimages in $S(V)^G$ of
$\{x_1,...,x_{n-3},\bar{w}_{n-2}^{p^f},a^{p^g},b^{p^h}\}$, $\{x'_1,...,x'_{n-3},x'_{n-2},x'_{n-1},x'_n\}$.
These $n$ preimages form a homogenous s.o.p. of $S(V)^G$ and by (\ref{eq2:14}):
\begin{equation*}\label{eq2:20}
\begin{split}
   \prod_{i=1}^{n}(\deg x'_i) & =
      \prod_{i=1}^{n-3}(\deg x_i)(\deg x'_{n-2})(\deg x'_{n-1})(\deg x'_n) \\
     &= \prod_{i=1}^{n-3}(\deg x_i) p^f p^g (\deg a)p^h(\deg b)
     = |T| (\deg a)(\deg b)\\
     & = |T||G/T| = |G|.
\end{split}
\end{equation*}
Therefore by Proposition \ref{P2-11.5} $S(V)^G=F[x'_1,...,x'_{n-3},x'_{n-2},x'_{n-1},x'_n]$ is a polynomial ring.$\Box$
\begin{Definition} \label{D2-16}
Let $A$ be an affine domain over a field $F$ and $\mathfrak{p} \in \spec A$. We say that $A_{\mathfrak{p}}$ is smooth if $A_{\mathfrak{p}}$ is geometrically regular, that is if $A_{\mathfrak{p}}\otimes _F K$ is regular for each finite field extension $K$ of $F$.
\end{Definition}
\begin{Remark}\label{R2-17} Note the following:
\begin{enumerate}
  \item By \cite[Theorem 28.7]{Matsumura} this is equivalent to $A_\mathfrak{\mathfrak{p}}$ being $\mathfrak{p}_\mathfrak{p}$-smooth over $F$;
  \item \cite[X, $\S6$, Definition 1, p. 75]{BOU-AC} uses the terminology "absolutely regular". It is slightly differently defined, but is an equivalent notion.
\end{enumerate}
\end{Remark}

\begin{Lemma}\label{L2-18}
Let $A$ be an affine domain over a field $F$. Then:
\begin{enumerate}
  \item $A_{\mathfrak{p}}$ is geometrically regular if and only if $\hat{A}_\mathfrak{p}$ is geometrically regular;
  \item Let $A \subset B$ be a finite extension, and $\mathfrak{p} \in \spec A$, $\mathfrak{n} \in \spec B$ with $\mathfrak{n} \cap A= \mathfrak{p}$. Suppose $B_{\mathfrak{n}}$ is \'{e}tale over $A_{\mathfrak{p}}$.
      Then $A_\mathfrak{p}$ is geometrically regular if and only if $B_{\mathfrak{p}}$ is geometrically regular.
\end{enumerate}
\end{Lemma}
\textbf{Proof:}
$(1)$ follows from \cite[$\S$28, Example 1(ii), p.215]{Matsumura}.

$(2)$ Let $K \supset F$ be a finite radical field extension of $F$.
$A_{\mathfrak{p}} \otimes _F K \subset B_{\mathfrak{n}} \otimes _F K$ is a faithfully flat extension \cite[p. 46, (2)]{Matsumura}.
Both rings are local, hence, by using \cite[Thm 23.7(i)]{Matsumura}, the regularity of $A_{\mathfrak{p}} \otimes _F K$ follows from the given regularity of $B_{\mathfrak{p}} \otimes _F K$.
This shows by Remark \ref{R2-17}(2) that $A_{\mathfrak{p}}$ is geometrically regular.

To show the opposite direction, suppose that $A_{\mathfrak{p}}$ is geometrically regular.
Recall that $\mathfrak{p}_{\mathfrak{p}}B_\mathfrak{n}=\mathfrak{n}_{\mathfrak{n}}$ and $B_{\mathfrak{n}}/\mathfrak{n}_{\mathfrak{n}}$ is a finite separable extension of $A_{\mathfrak{p}}/\mathfrak{p}_{\mathfrak{p}}$.
Consider the sequence $F \subset A_{\mathfrak{p}} \subset B_{\mathfrak{n}}$.
By \cite[Thm. 28.1 "transitivity"]{Matsumura} it suffices to show that $B_{\mathfrak{n}}$  is $\mathfrak{n}_{\mathfrak{n}}$-smooth over $A_{\mathfrak{p}}$. For this purpose we apply \cite[Thm. 28.10]{Matsumura} observing that
$B_{\mathfrak{n}}\otimes_{A_{\mathfrak{p}}} A_{\mathfrak{p}}/\mathfrak{p}_{\mathfrak{p}} =
B_{\mathfrak{n}}/\mathfrak{p}_{\mathfrak{p}}B_{\mathfrak{n}}=B_{\mathfrak{n}}/\mathfrak{n}_{\mathfrak{n}}$ is $0$-smooth over $A_{\mathfrak{p}}/\mathfrak{p}_{\mathfrak{p}}$, being a finite separable field extension of $A_{\mathfrak{p}}/\mathfrak{p}_{\mathfrak{p}}$.$\Box$

The next result is a complement to Lemma \ref{L2-4}.

\begin{Lemma} \label{L2-21}
Let $G \subset GL(V)$ be a finite group and $\dim _F V=n$.
Let $n \ge m \ge 1$ be an integer.
Consider the following properties:
\begin{enumerate}
\setcounter{enumi}{2}
  \item $S(V)^{G_U}$ is a polynomial ring for each subspace $U \subset V^*$ with $dim_FU \ge n-m$;
  \item $S(V)^G_{\mathfrak{p}}$ is geometrically regular for each $\mathfrak{p} \in \spec S(V)^G$ with $\height (\mathfrak{p})=m$.
  \end{enumerate}
Then we have the following implications: $(3) \Leftrightarrow (4)$.
\end{Lemma}
\textbf{Proof:}
Firstly we consider $(3) \Rightarrow (4)$.
Let $\mathfrak{\mathfrak{p}} \in \spec S(V)^G$, $\mathfrak{m} \in \spec S(V)$ with $\mathfrak{m} \cap S(V)^G=\mathfrak{p}$.
Set $\mathfrak{n}:=\mathfrak{m} \cap S(V)^{I_G(\mathfrak{m})}$, $\mathfrak{q}:=\mathfrak{m} \cap S(V)^{D_G(\mathfrak{m})}$.
Set $W:=\mathfrak{m} \cap V$ and $U:=W^{\perp}$.
Since $WS(V) \subseteq \mathfrak{m}$ is a prime ideal, we have that $\dim_F W = \height (WS(V)) \le \height (\mathfrak{m}) =\height (\mathfrak{p})=m$.
Consequently $\dim _F U=n - \dim_F W \ge n-m$.
Therefore by condition (3) $S(V)^{G_U}$ is a polynomial ring, implying that $S(V)^{G_U}$ is geometrically regular (every scalar field extension of a polynomial ring is a polynomial ring).
Since $G_U=I_G(\mathfrak{m})$ we have that $S(V)_{\mathfrak{n}}^{I_G(\mathfrak{m})}$ is geometrically regular, implying by Lemma \ref{L2-18}(2) and Theorem \ref{T2-2}(2) that $S(V)_{\mathfrak{q}}^{D_G(\mathfrak{m})}$ is geometrically regular. This implies by Lemma \ref{L2-18}(1) and Theorem \ref{T2-2}(1) that $S(V)_{\mathfrak{p}}^G$ is geometrically regular.

We now show the direction $(4) \Rightarrow (3)$.
By Lemma \ref{L2-4} we may only consider $U \subset V^*$ with $\dim _FU = n-m$.
We have to show that $S(V)^{G_U}$ is a polynomial ring.
Set $W:=U^{\perp}$, $\mathfrak{m}=WS(V)$.
Then $\mathfrak{m}$ is a prime ideal in $S(V)$ and $m=\dim _FW= \height (\mathfrak{m})$.
We also set $\mathfrak{n} := S(V)^{I_G(\mathfrak{m})}\cap \mathfrak{m}$, $\mathfrak{p}:= S(V)^G\cap \mathfrak{m}$, $\mathfrak{q}:=S(V)^{D_G(\mathfrak{m})} \cap \mathfrak{m}$.
We have $\height (\mathfrak{p}) = \height (\mathfrak{\mathfrak{m}})$ by "going down", and $I_G(\mathfrak{m})=G_U=\{g \in G|(g-1)(V) \subseteq W\}$.
By assumption $S(V)^G_{\mathfrak{p}}$ is geometrically regular and consequently by Theorem \ref{T2-2}(1) $\widehat{S(V)^G_{\mathfrak{p}}} \cong \widehat{S(V)^{D_G(\mathfrak{m})}_{\mathfrak{q}}}$ is geometrically regular.
Consequently $S(V)^{D_G(\mathfrak{m})}_{\mathfrak{q}}$ is geometrically regular and by Lemma \ref{L2-18}(2) so is $S(V)^{I_G(\mathfrak{m})}_{\mathfrak{n}} =S(V)^{G_U}_{\mathfrak{n}}$.
Consequently by \cite[X,\S6, Cor. 4, p. 77]{BOU-AC} $S(V)^{G_U}_{\mathfrak{n}} \otimes _{F} \bar{F} $ is  regular, where  $\bar{F}$ is the algebraic closure of $F$.

Let $V_{\bar{F}}:=V \otimes _{F} \bar{F}$, $W_{\bar{F}} =W \otimes _{F} {\bar{F}}$, $\mathfrak{m}_{\bar{F}} =\mathfrak{m} \otimes _F \bar{F}$.
Since $S(V_{\bar{F}})$ is integral over $S(V)^{I_G(\mathfrak{m})}$, there exists $\mathfrak{b} \in \spec S(V_{\bar{F}})$ such that
 $\mathfrak{b} \cap S(V)^{I_G(\mathfrak{m})} =\mathfrak{n}$ and $\mathfrak{b} \supseteq \mathfrak{m}_{\bar{F}}$.
Set $\mathfrak{a}=\mathfrak{b}\cap S(V_{\bar{F}})^{I_G(\mathfrak{m})}$, $H:=I_G(\mathfrak{m}) =G_U $.

Since $(g-1)(V_{\bar{F}}) \subseteq W_{\bar{F}} =V_{\bar{F}} \cap \mathfrak{m}_{\bar{F}} \subseteq V_{\bar{F}} \cap \mathfrak{b}$, $\forall g \in H$, we have
$(g-1)(S(V_{\bar{F}})) \subseteq \mathfrak{b}$, $\forall g \in H$.
Consequently $H=I_H(\mathfrak{b}) = D_H(\mathfrak{b})$ and $S(V_{\bar{F}})^H_{\mathfrak{a}}$ is regular being the localization of $S(V_{\bar{F}})^H_{\mathfrak{n}}=S(V)^H_{\mathfrak{n}}\otimes _F \bar{F}$.
Since $\bar{F}$ is perfect, we can conclude by \cite[Prop. 1.9, Example. 1.10]{K-Loci} that $S(V_{\bar{F}})^H=S(V_{\bar{F}})^{G_U}=S(V)^{G_U} \otimes _{F} \bar{F}$ is a polynomial ring, and therefore so is $S(V)^{G_U}$.$\Box$

\underline{\textbf{The proof of Theorem B:}}

Note that by Lemma \ref{L2-21} and Lemma \ref{L2-4}, $\dim$(non-smooth locus of $S(V)^G$)$\le 1$ is equivalent to: $S(V)^{G_U}$ is a polynomial ring for each $U \subset V^*$ with $\dim _FU=2$.
We have that $S(V)^G$ is Cohen-Macaulay and $S(V)^G_{\mathfrak{p}}$ is geometrically regular for each prime $\mathfrak{p}$ with $\height (\mathfrak{p})=n-2$.
Consequently by \cite[X, $\S6$, Cor. 4, p.77]{BOU-AC} $S(V)^G_{\mathfrak{p}} \otimes _F \bar{F}$ is regular, where $\bar{F}$ is the algebraic closure of $F$.
Set $V_{\bar{F}}
= V \otimes _F \bar{F}$, then $S(V_{\bar{F}})_{\mathfrak{p}}^G
= S(V)^G_{\mathfrak{p}}\otimes_{F} \bar{F}$ is regular.

Let $\mathfrak{q} \in \spec S(V_{\bar{F}})^G$ with $\height (\mathfrak{q})=n-2$.
Since $S(V_{\bar{F}})^G$  is integral over $S(V)^G$, it follows that $\mathfrak{q} \cap S(V)^G$ is a height $(n-2)$ prime ideal in $S(V)^G$.
Consequently by the previous paragraph, $S(V_{\bar{F}})^G _{\mathfrak{q} \cap S(V)^G} =S(V)^G_{\mathfrak{q} \cap S(V)^G} \otimes _{F} \bar{F}$ is regular. Therefore $S(V_{\bar{F}})^G_\mathfrak{q}$, being a localization of the latter ring, is also regular.

To sum up, $S(V_{\bar{F}})^G$ satisfies Serre's $R_{n-2}$ condition, $\bar{F}$ is perfect and\\ $S(V_{\bar{F}})^G = S(V)^G \otimes_{F}\bar{F}$ is Cohen-Macaulay.
So by Lemma \ref{L2-4} all the requirements of Theorem \ref {T2-12} are in place and $S(V_{\bar{F}})^G$ is therefore a polynomial ring, implying that $S(V)^G$ is a polynomial ring.$\Box$

In the next example we exhibit a $p$-group $T \subset GL(V)$, $\dim _F V=4$, with the additional properties $|T|=p^3$, $T$ is elementary abelian and is generated by transvections.
However, $S(V)^T$ is not a polynomial ring, $S(V)^T$ is Gorenstein and $T=T_U$, for some subspace $U \subset V^*$ with $\dim U=2$.

\begin{Example} \label{EX2-1}
$T=<\sigma, \tau _1, \tau _2> \subset GL_4(F)$ and $F \supset \mathbb{F}_p$ where:
$$\sigma=
 \begin{bsmallmatrix}
    1 & 0 & 1 & 0 \\
    0 & 1 & 0 & 0 \\
    0 & 0 & 1 & 0 \\
    0 & 0 & 0 & 1
  \end{bsmallmatrix}, \text{ }
\tau _1=
 \begin{bsmallmatrix}
    1 & 0 & 0 & 1 \\
    0 & 1 & 0 & 1 \\
    0 & 0 & 1 & 0 \\
    0 & 0 & 0 & 1
  \end{bsmallmatrix}, \text{ }
\tau _2=
 \begin{bsmallmatrix}
    1 & 0 & 0 & 0 \\
    0 & 1 & 0 & \xi \\
    0 & 0 & 1 & 0 \\
    0 & 0 & 0 & 1
  \end{bsmallmatrix}, \text{ }
 \xi \in F -\mathbb{F}_p. $$
The matrices are given with respect to the basis $\{v, w_3, w_2, w_1 \}$.
We read the action of the matrices along the rows, so e.g. $\sigma (v)=v+w_2$, $\tau _1(v)=v+w_1$, $\tau _1(w_3)=w_3+w_1$, $\tau _2(w_3)=w_3+\xi w_1$, etc.

Set $W:=Fw_3+Fw_2+Fw_1$. We have:
$V^{\sigma}=Fw_3+Fw_2+Fw_1$, $V^{\tau _1}=F(v-w_3)+Fw_2+Fw_1$, $V^{\tau _2}=Fv+Fw_2+Fw_1$.
So $V^T=Fw_2+Fw_1$.
Clearly $|<\sigma , \tau _2>|=p^2$ and $S(V)^{<\sigma , \tau _2>} = F[w_1, w_2,w^p_3-(\xi w_1)^{p-1}w_3, v^p-w_2^{p-1}v]$.

Also $S(W)^T=S(W)^{<\tau _1, \tau _2>}$, so by Proposition \ref{P2-6.5} $S(W)^T=F[w_1,w_2,a]$, with $\deg a =|<\tau _1, \tau _2>|=p^2$.

Suppose $S(V)^T$ is a polynomial ring.
So $S(V)^T=F[w_1,w_2,x,y]$, and consequently $(\deg x )(\deg y) = |T|=p^3$ (this is so since $V^T=Fw_2+Fw_1$).
Hence by Proposition \ref{P2-11.5}, $\deg x=p$, $\deg y=p^2$.
At least one of $\{x,y\}$ must be in $S(V)^T-S(W)$, otherwise $\{w_1,w_2,x,y\}$ are algebraically dependent over $F$.

Suppose $x \in S(V)^T-S(W)$, then $S(V)^T \subset S(V)^{<\sigma, \tau _2>}$ and $x \in S(V)^{<\sigma, \tau _2>}$.
Therefore since $\deg x=p$, we have $x=\beta (w_1,w_2)+\lambda [w_3^p - (\xi w_1)^{p-1}w_3]+\mu [v^p-w_2^{p-1}v]$, $\mu \neq 0$, $\lambda, \mu \in F$, $\beta (w_1, w_2) \in F[w_1,w_2]$, $\deg \beta (w_1,w_2)=p$. Hence:
\begin{equation*}\label{eq2:20}
\begin{split}
0 & = (\tau _1 -1)(x) \\
  & =  0+\lambda [(\tau _1-1)(w_3)^{p}-(\xi w_1)^{p-1}(\tau _1 -1)(w_3)]+ \mu [(\tau _1 -1)(v)^p - (w_2)^{p-1}(\tau _1 - 1)(v)] \\
  & = \lambda [w_1^p -(\xi w_1)^{p-1}w_1] + \mu [w_1^p-w_2^{p-1}w_1] = \lambda (1 - \xi^{p-1})w_1^p + \mu [w_1^p-w_2^{p-1}w_1],\end{split}
\end{equation*}
a contradiction, since $\mu \neq 0$ and $\{w_1,w_2\}$ are algebraically independent.

Therefore $x \in S(W)^T$.
By the argument above $S(W)^T=F[w_1,w_2,a]$, where $\deg a = p^2$.
Since $\deg x = p$ we have $x \in F[w_1,w_2]$, which is another contradiction since $\{w_1, w_2,x, y \}$ is a minimal generating set of $S(V)^T$,
so $x$ must be algebraically independent from $\{w_1,w_2\}$.

Finally we show that $S(V)^T$ is a Gorenstein ring.
Let $\mathfrak{m} := (w_1,w_2)S(V), \mathfrak{n}:= \mathfrak{m} \cap S(V)^{I_T(\mathfrak{m})}$, $\mathfrak{p}:=\mathfrak{m} \cap S(V)^T$.
We have that $V \cap \mathfrak{m} = Fw_2+Fw_1$ and consequently $T_U=I_T(\mathfrak{m})$, where $U:=(Fw_2+Fw_1)^{\perp} \subset V^*$.
Since $S(V)$ is integral over $S(V)^{T_U}$, it follows that $\height (\mathfrak{n}) = \height (\mathfrak{m})=2$.
Moreover since $S(V)^T$ is normal, it follows that $S(V)^T_{\mathfrak{p}}$ is Cohen-Macaulay, implying by \cite[Theorem 1.9]{K-Loci} that $S(V)^{T_U}$ is Cohen-Macaulay.
However by inspection $T_U=\{g \in T|(g -1)(V) \subseteq Fw_2+Fw_1\}=T$, and hence $S(V)^T$ is Cohen-Macaulay.
Since $T$ is generated by pseudo-reflections, it is a U.F.D, implying by the Cohen-Macaulay property, that it is Gorenstein.
\end{Example}

The nature of $S(V^*)^T$ is surprisingly different.
\begin{Example} \label{EX2-2}
Let $T \subset GL(V)$ be as in Example \ref{EX2-1}.
Then $S(V^{*})^T$ is a polynomial ring, and $T=T_{Fw_1+Fw_2}$.
\end{Example}

\textbf{Proof:} Recall that the matrices of $\{\sigma, \tau _1, \tau _2\}$ in $GL(V)$ where given in Example \ref{EX2-1} with respect to the ordered basis $\{v,w_3,w_2,w_1\}$.

Let $\{g, f_3,f_2,f_1\} \subset V^*$ be the dual basis of $\{w_1,w_2,w_3,v\}$, that is:\\
 $f_1(v)=1$, $f_1(w_3)=f_1(w_2)=f_1(w_1)=0$; $f_2(w_3)=1$, $f_2(v)=f(w_2)=f_2(w_1)=0$; $f_3(w_2)=1$, $f_3(v)=f_3(w_3)=f_3(w_1)=0$; $g(w_1)=1$, $g(v)=g(w_3)=g(w_2)=0$.

With respect to this ordered basis the matrices of $\sigma$, $\tau _1$, $\tau _2$ are given by:
$$\sigma=
 \begin{bsmallmatrix}
    1 & 0 & 0 &  0 \\
    0 & 1 & 0 & -1 \\
    0 & 0 & 1 &  0 \\
    0 & 0 & 0 &  1
  \end{bsmallmatrix}, \text{ }
\tau _1=
 \begin{bsmallmatrix}
    1 & 0 & -1 & -1 \\
    0 & 1 &  0 &  0 \\
    0 & 0 &  1 &  0 \\
    0 & 0 &  0 &  1
  \end{bsmallmatrix}, \text{ }
\tau _2=
 \begin{bsmallmatrix}
    1 & 0 & -\xi & 0 \\
    0 & 1 &  0   & 0 \\
    0 & 0 &  1   & 0 \\
    0 & 0 &  0   & 1
  \end{bsmallmatrix}.$$
Set $a_1:=g^p-(\tau _1 -1)(g)^{p-1}g=g^p-(-f_2 - f_1)^{p-1}g$.
Clearly $\tau _1(a_1)=a$.
Set $a = a_1^p-(\tau _2-1)(a_1)^{p-1}a_1=a_1^p -[(-\xi f_2)^p-(-f_2-f_1)^{p-1}(-\xi f_2)]^{p-1}a_1$.
It is easily seen to be a $p$-polynomial of degree $p^2$ in $g$, and $(\sigma -1)(a)=0$ holds,
since $(\sigma - 1)(g)= (\sigma -1)(f_2)=(\sigma - 1)(f_1) =0$.
Also $\tau _1(a)=a$, $\tau _2(a)=a$.

Therefore $\{f_1,f_2,f_3,a\}$ are homogenous elements in $S(V^*)^{<\tau _1, \tau _2>}$, as well as algebraically independent over $F$ and
$(\deg f_1)(\deg f_2)(\deg f_3)(\deg a) =\deg a =|<\tau _1, \tau _2>|=p^2$.
Hence $S(V^*)^{<\tau _1, \tau _2>} =S(U_1)[a] =F[f_1, f_2, f_3, a]$, where $U_1:=Ff_1+Ff_2+Ff_3$.
Consequently $S(V^*)^{<\tau _1, \tau _2, \sigma>}=(S(V^*)^{<\tau _1, \tau _2>})^{<\sigma>} =(S(U_1)[a])^{<\sigma>}=S(U_1)^{<\sigma >}[a]=F[f_1,f_2,f_3^p-(-f_1)^{p-1}f_3,a]$, which is a polynomial ring.$\Box$

\section{\bf Theorem A}\label{SEC-3}

Recall that we adhere to the convention of acting with matrices on the right, so the action is read along rows.
\begin{Proposition}\label{P3-1}
Let $dim_FV=3$, $p=charF >0$ and $T \subset GL(V)$ a finite transvection group, $T=<\sigma, \tau _1, ... , \tau _m, \epsilon _1,... ,\epsilon_n>$, where with respect to a fixed basis $\{v, w_2, w_1\}$, we have:
\begin{enumerate}
\item $\tau _i :=
\begin{bsmallmatrix}
1 & 0 & \tau _{13}^{(i)} \\
0 & 1 & \tau _{23}^{(i)} \\
0 & 0 & 1
\end{bsmallmatrix}$,
$\tau_{23}^{(i)}\neq 0$, $i=1,...,m$,
\item $\sigma :=
\begin{bsmallmatrix}
1 & \sigma _{12} & \sigma _{13} \\
0 & 1 & 0 \\
0 & 0 & 1
\end{bsmallmatrix}$,
$\sigma_{12} \neq 0$,
\item $\epsilon _i:=
\begin{bsmallmatrix}
1 & 0 &\epsilon _{13}^{(i)} \\
0 & 1 & 0 \\
0 & 0 & 1
\end{bsmallmatrix}$, $\epsilon _{13}^{(i)} \neq 0$, $i=1,...,n$.
\end{enumerate}
Suppose further:
\begin{enumerate}\addtocounter{enumi}{+3}
\item $(H,H) \subseteq <\epsilon _1, ... ,\epsilon _n >$, where $H:=<\sigma ,\tau _1, ... , \tau _m>$,
\item $\{\epsilon _1,... ,\epsilon _n\}$ is a minimal generating set of $<\epsilon_1, ... , \epsilon _n>$, equivalently $\{ \epsilon _{13}^{(i)} |i=1,...,n\}$ is $\mathbb{F}_p$-linearly independent.
\item $\{\tau _{1|W}, ... , \tau _{m|W}\}$ is a minimal generating set for the subgroup $<\tau _{1|W}, .... , \tau _{m|W} >$, equivalently $\{\tau _{23}^{(i)}|i=1, ... , m\}$ are $\mathbb{F} _p$-linearly independent.
\end{enumerate}
Then $S(V)^T$ is a polynomial ring.
\end{Proposition}

\textbf{Proof:} Set $W:=Fw_1+Fw_2$.
The plan of the proof is as follows.
We aim at showing that $S(V)^T=F[w_1,y,X]$, where $S(W)^T=F[w_1,y]$, $\deg y=|T_{|W}|$.
This crucially uses the fact that $m \le n$, where $p^m=|T_{|W}|$ and $p^n=|<\epsilon_1,...,\epsilon_n>|$.
Also $S(V)^{<\sigma, \epsilon_1,...,\epsilon_n>}=F[w_1,w_2,a]$, $\deg a =p^{n+1}$ and $X$ is obtained from the element $a$ by a sequence of $m$ adjustments.
Since $(\deg w_1)(\deg y)(\deg X)= 1\cdot p^m \cdot p^{n+1} =|T|$, the result will follow from Proposition \ref{P2-11.5}.

Now
$\sigma \tau _i \sigma ^{-1}\tau _i^{-1}=
\begin{bsmallmatrix}
1 & 0 &\sigma _{12}\tau^{(i)} _{23} \\
0 & 1 & 0 \\
0 & 0 & 1
\end{bsmallmatrix}$, $i=1,... , m$, and $\sigma _{12} \neq 0$, so by (6) $\{\sigma_{12} \tau_{23}^{(i)}|i=1,..., m\}$ are $\mathbb{F}_p$-linearly independent, implying by (4), (5) that $m \le n$.

Set $z_0:=v$, $z_{i+1}:=z_{i}^p-(\epsilon _{i+1}-1)(z_i)^{p-1}z_i$, $i=0,... ,n-1$.
Then by Proposition \ref{P2-6.5} $z:=z_n$ is of degree $p^n$ and $S(V)^{<\epsilon _1, ..., \epsilon _n>}=S(W)[z]=F[w_2,w_1,z]$.
Set $a:=z^p-(\sigma -1)(z)^{p-1}z$.
Since $\deg(z)=p^n$ we have $\deg(a)=p^{n+1}$.

Similarly let $y_{i+1}:=y_i^p-(\tau _{i+1}-1)(y_i)^{p-1}y_i$, $i=0,...,m-1$, $y_0=w_2$.
Clearly $\deg(y_i)=p^i$.
Then by Proposition \ref{P2-6.5} $S(W)^{<\tau _1,... ,\tau_i>}=F[w_1,y_i]$, $i=1, ... ,m$.
Set $y:=y_m$, so $\deg(y)=p^m$.

Moreover:
\begin{enumerate}\addtocounter{enumi}{+6}
\item $z_i=v^{p^i}+\sum_{j=0}^{i-1}\eta_j v^{p^j}w_1^{p^i-p^j}$, $\eta _j\in F$, $i=0,... ,n$;
\item $(\epsilon _j -1)(z_i)=\mu _{j,i}w_1^{p^i}$, $\mu_{j,i}\in F_i$, $j=1, ... ,n$, $i=0,... ,n$;
\item $(\tau _l -1)(y_j)=\theta _{l,j}w_1^{p^j}$, $\theta _{l,j} \in F$, $l=1, ... , m$, $j=0, ...,m$;
\item $(\tau _l-1)(z_i)=\alpha_{l,i}w_1^{p^i}$, $\alpha _{l,i} \in F$, $l=1,...,m$, $i=0,... ,n$.
\end{enumerate}


We shall next show that $(\tau _l-1)((\sigma - 1)(z))=0$, $l=1,... ,m$.

Indeed $(\tau _l \sigma)(z)=(\tau _l\sigma \tau_l ^{-1}\sigma ^{-1})(\sigma \tau _l)(z)=(\sigma \tau_l)((\tau _l\sigma \tau _l^{-1}\sigma ^{-1})(z))=(\sigma \tau _l)(z)$, where the second equality follows since $\tau _l\sigma \tau _l^{-1} \sigma ^{-1}$ is in  $<\epsilon _1, ... ,\epsilon _n> \subset Z(T)$ and the last equality holds since $\epsilon _i(z)=z$ for each $i$.
Hence $(\tau _l -1)((\sigma -1)(z))=(\sigma - 1)((\tau _l -1)(z))$.
Thus by (10), $(\sigma -1)((\tau _l-1)(z))=(\sigma - 1)(\alpha _{l,n}w_1^{p^n})=0$, for $l=1, ... ,m$.
Consequently $(\tau _l-1)((\sigma -1)(z))=0$, $l=1, ... ,m$, as claimed.

Since $(\sigma -1)(z)\in F[w_1,w_2]=S(W)$, it follows that $(\sigma -1)(z) \in S(W)^{<\tau_1,...,\tau _m>}=F[w_1,y]$.

Assume, by induction, the existence of $b_i$ of the following form:\\
$b_i=\xi _iw_1^{p^{n+1}-p^{i-1}}y_{i-1}-\eta _i[(\sigma -1)(z)]^{p-1}w_1^{p^n-p^{i-1}}y_{i-1}$, $\xi _i,\eta _i \in F$, and satisfying:
\begin{enumerate}\addtocounter{enumi}{+10}
\item $\tau _i(a-b_1-\cdots -b_i)=a-b_1-\cdots -b_i$.
\end{enumerate}

Since $(\sigma -1)(z)\in S(W)^{<\tau_1, ... ,\tau _m>}=F[w_1, y_m]$, it follows that $b_i\in F[w_1, y_{i-1}]=S(W)^{<\tau_1, ... ,\tau _{i-1}>}$, $i=2,... , m$.

We shall next construct $b_{i+1}\in F[w_1,y_i]$, in the above form and satisfying:\\
$\tau _{i+1}(a-b_1- \cdots -b_{i+1})=a-b_1- \cdots -b_{i+1}$.
We have by (9):
\begin{equation*}\label{eq3:1}
\begin{split}
   (\tau_{i+1}-1)(b_j) & =\xi _jw_1^{p^{n+1}-p^{j-1}}(\tau _{i+1}-1)(y_{j-1})-\eta _j [(\sigma - 1)(z)]^{p-1}w_1^{p^n-p^{j-1}}(\tau _{i+1}-1)(y_{j-1}))\\
     & = \xi _j\theta_{i+1,j-1}w_1^{p^{n+1}} - \eta _j\theta_{i+1,j-1}[(\sigma -1)(z)]^{p-1}w_1^{p^n}, \text{ } j=1,.... ,i.
\end{split}
\end{equation*}
Therefore by using (10):
\begin{enumerate}\addtocounter{enumi}{+11}
\item $(\tau _{i+1}-1)(a-b_1-\cdots -b_i) $\\
 $ = (\tau _{i+1}-1)(a)-\sum_{j=1}^{i}(\tau_{i+1}-1)(b_j)$\\
 $ =  [(\tau_{i+1}-1)(z)^p-(\sigma -1)(z)^{p-1}(\tau _{i+1}-1)(z)]$\\
 $ \indent - \sum_{j=1}^{i}\{\xi _j\theta_{i+1,j-1}w_1^{p^{n+1}}-\eta _j \theta_{i+1,j-1}[(\sigma - 1)(z)]^{p-1}w_1^{p^n}\}$\\
 $ = (\alpha_{i+1,n}^p -\sum_{j=1}^{i}\xi _j\theta _{i+1,j-1})w_1^{p^{n+1}}-[(\sigma -1)(z)]^{p-1}(\alpha _{i+1,n} -\sum_{j=1}^{i}\eta _j\theta _{i+1,j-1})w_1^{p^n}$.
\end{enumerate}
Recall that by (9) $(\tau _{i+1}-1)(y_i)=\theta_{i+1,i}w_1^{p^i}$.
If $\theta _{i+1,i}=0$ then $S(W)^{<\tau_1,... ,\tau_i>} = F[w_1, y_i] \subseteq S(W)^{<\tau_1,... ,\tau_{i+1}>}$ in violation of (6).

Set $\xi_{i+1}:=\frac{\alpha _{i+1,n}^p-\sum_{j=1}^{i}\xi_j \theta_{i+1,j-1}}{\theta_{i+1,i}}$, $\eta_{i+1}:= \frac{\alpha _{i+1,n}-\sum_{j=1}^{i}\eta_j \theta_{i+1,j-1}}{\theta_{i+1,i}}$ and:\\
$b_{i+1}:=\xi_{i+1}w_1^{p^{n+1}-p^i}y_i - \eta_{i+1}[(\sigma - 1)(z)]^{p-1}y_iw_1^{p^n-p^i}.$\\
So it follows by (12) that $(\tau_{i+1}-1)(b_{i+1})=(\tau_{i+1}-1)(a-b_1-\cdots -b_i)$, and therefore $(\tau_{i+1}-1)(a-b_1- \cdots - b_i-b_{i+1})=0$.
So the inductive step is verified.
Note that $b_{i+1}=0$ is also possible.

Consequently $a-b_1-\cdots -b_i \in S(V)^{<\sigma, \epsilon_1, ..., \epsilon _n, \tau _1,... ,\tau_i>}$ for each $i=1,... ,m$.
Also $\deg(a-b_1-\cdots -b_i)=p^{n+1}$, $y_i \in S(V)^{<\sigma, \epsilon_1, ..., \epsilon _n, \tau _1,... ,\tau_i>}$ and $\deg(y_i)=p^i$.
Now $\sigma_{12}\neq 0$, so (5) implies that $|<\sigma, \epsilon_1,... ,\epsilon _n>|=p^{n+1}$.

We shall next observe that $<\tau_1,... ,\tau_m>\cap <\sigma, \epsilon_1,... ,\epsilon _n> =1$.
Indeed if $\tau_1^{e_1}\cdots \tau_m^{e_m} \in <\sigma, \epsilon_1,... ,\epsilon _n>$, with $ 0 \le e_i \le p-1$, then this element acts trivially on $S(W)$.
Hence $\tau_{1|W}^{e_1} \cdots \tau_{m|W}^{e_m}=0$, implying by (6) that $e_1=\cdots =e_m=0$.
Consequently $<\sigma, \epsilon_1,... ,\epsilon _n, \tau _1,... ,\tau_i>=p^{n+1}p^i$.
Since $\{w_1$, $y_i$, $a-b_1-\cdots -b_i\}$ are $F$-algebraically independent, it follows by Proposition \ref{P2-11.5} that $F[w_1,y_i,a-b_1-\cdots -b_i]=S(V)^{<\sigma, \epsilon_1,... ,\epsilon _n, \tau _1,... ,\tau_i>}$, for $i=1,... ,m$. Therefore $S(V)^T=F[w_1,y,X]$, where $X:=a-b_1 -\cdots -b_m$. $\Box$

The next result extends Proposition \ref{P3-1}. The basic idea is as in Proposition \ref{P2-7}.

\begin{Proposition}\label{P3-2}
Let $dim_FV=3$, $p=charF >0$ and $T \subset GL(V)$ the finite transvection group, $T=<\sigma_r, ... ,\sigma_1, \tau _1, ... , \tau _m, \epsilon _1,... ,\epsilon_n>$, where with respect to a fixed basis $\{v, w_2, w_1\}$, we have:
\begin{enumerate}
\item $\tau _i :=
\begin{bsmallmatrix}
 1 & 0 &\tau _{13}^{(i)} \\
 0 & 1 & \tau _{23}^{(i)} \\
 0 &0 &1 \end{bsmallmatrix}$,
$\tau_{23}^{(i)}\neq 0$, $i=1, ... ,m$,
\item $\sigma_i :=
\begin{bsmallmatrix}
 1 & \sigma _{12}^{(i)} &\sigma _{13}^{(i)} \\
 0 & 1 & 0 \\
 0 &0 & 1
 \end{bsmallmatrix}$,
 $\sigma_{12}^{(i)} \neq 0$, $i=1, ... ,r$,
\item $\epsilon _i:=
\begin{bsmallmatrix}
 1 & 0 &\epsilon _{13}^{(i)} \\
 0 & 1 & 0 \\
 0 & 0 & 1
\end{bsmallmatrix}$, $\epsilon _{13}^{(i)} \neq 0$, $i=1,... ,r$.
\end{enumerate}
Suppose further:
\begin{enumerate}\addtocounter{enumi}{+3}
\item $(H,H) \subseteq <\epsilon _1, ... ,\epsilon _n >$, where $H:=<\sigma_r, ... , \sigma_1, \tau _1, ... , \tau _m>$,
\item $\{\epsilon _1,... ,\epsilon _n\}$ is a minimal generating set of $<\epsilon_1, ... , \epsilon _n>$.
\end{enumerate}
Then $S(V)^T$ is a polynomial ring.
\end{Proposition}
\textbf{Proof:}
We may also assume:
\begin{enumerate}\addtocounter{enumi}{+5}
\item $\{\tau _{1|W}, ... , \tau _{m|W}\}$ is a minimal generating set of $<\tau _{1|W}, .... , \tau _{m|W} >$.
Equivalently $\{\tau_{23}^{(i)}|i=1, ... ,m\}$ is $\mathbb{F}_p$-linearly independent.
\end{enumerate}

Indeed say $\{\tau _{1|W}, ... , \tau _{r|W}\}$, $r <m$ is a minimal generating set for this subgroup. So for $i>r$ we have $ \tau _{i|W} \in <\tau _{1|W}, ... , \tau _{r|W}>$.
Therefore $\tau _{i|W}=\tau _{1|W}^{e_1}\cdots \tau _{r|W}^{e_r}$, $0 \le e_j \le p-1$, $1 \le j \le r$ implying that
$\delta _i=\tau _i \cdot \tau_1^{-e_1} \cdots \tau _r^{-e_r} \in
\begin{bsmallmatrix}
1 & 0 & F \\
0 & 1 & 0 \\
0 & 0 & 1
\end{bsmallmatrix}$.
Consequently we can replace $\tau _i$ by adding $\delta _i$ to $\{\epsilon _1, ... ,\epsilon _n\}$ in the generating set.

Similarly we may assume:
\begin{enumerate}\addtocounter{enumi}{+6}
\item $\{(\sigma_1) _{|V/Fw_1},... ,(\sigma_r) _{|V/Fw_1}\}$ is a minimal generating set of $<(\sigma_1) _{|V/Fw_1},... ,(\sigma_r) _{|V/Fw_1}>$.
Equivalently $\{\sigma_{12}^{(i)}|i=1,... ,r\}$ are $\mathbb{F}_p$-linearly independent.
\end{enumerate}

As a result of (4), (5), (6) and (7) we have $|<\sigma _j,... ,\sigma_1, \tau_1,... , \tau_m,\epsilon_1,... ,\epsilon_n>|=p^j\cdot p^m \cdot p^n$, for $j=1,... ,r$.

By Proposition \ref{P3-1} we have $S(V)^{<\sigma_1, \tau _1, ... , \tau _m, \epsilon _1,... ,\epsilon_n>}=F[w_1, y, X_1]$, where $y \in S(W)^T$,
$X_1=a-b_1-\cdots -b_m$ of Proposition \ref{P3-1}, $\deg(X_1)=p^{n+1}$.

Define by induction $X_{j+1}:=X_j^p-(\sigma_{j+1}-1)(X_j)^{p-1}X_j$, $j=1,... ,r-1$.
Assume by induction that $F[w_1, y,X_i]=S(V)^{<\sigma_i, ... , \sigma_1, \tau _1, ... , \tau _m, \epsilon _1,... ,\epsilon_n>}$.

We shall now verify the equality $F[w_1,y,X_{i+1}]=S(V)^{<\sigma_{i+1}, ... , \sigma_1, \tau _1, ... , \tau _m, \epsilon _1,... ,\epsilon_n>}$.
It follows from (4) that $(T,T) \subseteq <\epsilon_1,... ,\epsilon_n>\subseteq Z(T)$ and therefore $<\sigma_i, ... , \sigma_1, \tau _1, ... , \tau _m, \epsilon _1,... ,\epsilon_n>$ is normal in $T$.
Consequently $S(V)^{<\sigma_i, ... , \sigma_1, \tau _1, ... , \tau _m, \epsilon _1,... ,\epsilon_n>}$ is $\sigma _{i+1}$-stable, so $\sigma_{i+1}(X_i) \in F[w_1,y,X_i]$.
Since $\deg(\sigma_{i+1}(X_i))=\deg(X_i)=p^{(n+1)+i-1}>p^m=\deg(y)$, it follows that $\sigma_{i+1}(X_i)=\alpha X_i +t$, $t\in F[w_1,y]$, $\alpha \in F$.
Since $\sigma_{i+1}^p=1$ and $F[w_1,y]\subseteq F[w_1,w_2]$ is fixed by $\sigma_{i+1}$, it follows by iterations that $\alpha ^p=1$, so $\alpha =1$.
Therefore $(\sigma_{i+1}-1)((\sigma_{i+1}-1)(X_i))=(\sigma_{i+1}-1)(t)=0$.
Hence:
\begin{equation*}
\begin{split}
   (\sigma_{i+1}-1)(X_{i+1}) & = (\sigma_{i+1}-1)[(X_i^p-(\sigma_{i+1}-1)(X_i)^{p-1}X_i] = (\sigma _{i+1}-1)[X_i^p-t^{p-1}X_i] \\
     & =  [(\sigma_{i+1}-1)(X_i)]^p-t^{p-1}(\sigma_{i+1}-1)(X_i)= t^p-t^p=0.
\end{split}
\end{equation*}
Therefore $\{X_{i+1},y,w_1\}\subset S(V)^{<\sigma_i+1, ... , \sigma_1, \tau _1, ... , \tau _m, \epsilon _1,... ,\epsilon_n>}$ are algebraically independent and $(\deg(X_{i+1}))(\deg(y))(\deg(w_1))=p^{(n+1)+i}\cdot p^m \cdot 1=|<\sigma_{i+1}, ... , \sigma_1, \tau _1, ... , \tau _m, \epsilon _1,... ,\epsilon_n>|$.
So again by Proposition \ref{P2-11.5} $F[w_1,y,X_{i+1}]=S(V)^{<\sigma_{i+1}, ... , \sigma_1, \tau _1, ... , \tau _m, \epsilon _1,... ,\epsilon_n>}$.
Hence $S(V)^T=F[w_1,y,X]$, where $X:=X_n$.$\Box$

We finally arrive at the main result of this section.

\begin{Theorem}\label{T3-3}
Suppose $T \subset GL(V)$ is a finite transvection group, $dim_FV=3$ and $Fw_1 \subset W = Fw_2+Fw_1 \subset V$ is a flag of $T$-modules.
Then $S(V)^T$ is a polynomial ring.
Moreover $S(V)^T=F[w_1,y,X]$, with $y \in S(W)$, if $T \not\subset
\begin{bsmallmatrix}
 1 & 0 & F \\
 0 & 1 & F \\
 0 & 0 & 1
\end{bsmallmatrix}$,
$\begin{bsmallmatrix}
 1 & F & F \\
 0 & 1 & 0 \\
 0 & 0 & 1
\end{bsmallmatrix}$.
\end{Theorem}
\textbf{Proof:}
Let $\theta \in T$ be a transvection.
Then $\theta(w_1)=w_1$.
Also $\theta_{|W}$ is a transvection on $W$ (or $1_W$) implying that $\theta_{|W}\in
\begin{bsmallmatrix}
 1 & F \\
 0 & 1
\end{bsmallmatrix}$ (with respect to the basis $\{w_2,w_1\}$).
Since $det \theta=1$ it follows that $\theta \in
\begin{bsmallmatrix}
 1 & F & F \\
 0 & 1 & F \\
 0 & 0 & 1
\end{bsmallmatrix}$, with respect to the basis $\{v,w_2,w_1\}$ where $v \in V - W$.
Since $rank(\theta -1_V)=1$, it follows that either $\theta\in
\begin{bsmallmatrix}
 1 & F & F \\
 0 & 1 & 0 \\
 0 & 0 & 1
\end{bsmallmatrix}$ or $\theta \in
\begin{bsmallmatrix}
 1 & 0 & F \\
 0 & 1 & F \\
 0 & 0 & 1
\end{bsmallmatrix}$.
Therefore each transvection generator of $T$ has the form of either (1), (2) or (3) of Proposition \ref{P3-2}.

If $\sigma _r=\cdots =\sigma _1=1$, then by Proposition \ref{P2-7} $S(V)^T$ is a polynomial ring.
Similarly if $\tau _1=\cdots =\tau_m=1$, then by Proposition \ref{P2-6.5} $S(V)^T$ is again a polynomial ring.
So we may assume that $T$ is generated by transvections as described in Proposition \ref{P3-2}.

Keeping the notation of Proposition \ref{P3-2}, then since $(H,H)$ is a finite subgroup of $(T,T) \subseteq
\begin{bsmallmatrix}
 1 & 0 & F \\
 0 & 1 & 0 \\
 0 & 0 &1
\end{bsmallmatrix}$
it follows that the generators of $(H,H)$ can be added to $<\epsilon_1,..., \epsilon_n>$, if they are not already there, so condition (4) of Proposition \ref{P3-2} can be arranged.
Finally condition (5) of Proposition \ref{P3-2} can also be arranged by dropping any redundant generator.
Thus all the requirements of Proposition \ref{P3-2} are in place, and $S(V)^T$ is therefore a polynomial ring.$\Box$

The next result is a direct consequence.

\begin{Theorem} \label{T2-4}
Let $T \subset GL(V)$ be a finite $p$-group which is generated by transvections, where $p= \Char F>0$ and $\dim _F V=3$.
Then $S(V)^T$ is a polynomial ring.
\end{Theorem}
\textbf{Proof:}
By a standard result for $p$-groups, there is a flag of $T$-modules $Fw_1\subset W =Fw_2+Fw_1 \subset V$.
Therefore Theorem \ref{T3-3} is applicable.$\Box$

\section{\bf Appendix}\label{Appendix}

In order to better explain the proof of Theorem \ref{T2-12} we present below a detailed sketch of the special case $\dim_FV=4$.
It avoids the somewhat complicated induction appearing in the proof, resulting in a much shorter proof where most crucial ingredients are still presented.

\underline{The proof of Theorem \ref{T2-12} when $\dim_FV=4$, $F=F^p$}

We assume that $S(V)^G$ is Cohen-Macaulay and satisfies Serre's $R_2$ condition.
$G$ is generated by transvections by Lemma \ref{L2-5}.
By a standard result there is a non-zero element $v_1 \in V$ with $g(v_1)=v_1$, for each $g \in G$.

Let $H_1:=\{g \in G|(g-1)(V)\subseteq Fv_1\}$, $V_2:=V/Fv_1$.
Then since $H_1=\ker (G \rightarrow G_{|V_2})$, $H_1$ is normal in $G$ and $G/H_1:=G_2$ acts faithfully and linearly on $V_2$.
It is a consequence of Proposition \ref{P2-7} that $S(V)^{H_1}=F[v_1,x_2,x_3,x_4]$ is a polynomial ring, where $\{x_2,x_3,x_4\}$ are homogenous elements, $\deg x_i=p^{f_i}$, and $x_i- u_i^{p^{f_i}} \in v_1S(V)$, $u_i \in V$, $i=2,3,4$.
Consequently by Proposition \ref{P2-11.5} $|H_1|=p^{f_2}\cdot p^{f_2}\cdot p^{f_4}$.

Since $S(V_2) \supset F[\bar{u}_2^{p^{f_2}},\bar{u}_3^{p^{f_3}},\bar{u}_4^{p^{f_4}}]= S(V)^{H_1}/v_1S(V)^{H_1}$ and the latter is a polynomial ring in $3$ variables, it follows that $\{\bar{u}_2,\bar{u}_3,\bar{u}_4\}$ is a basis of $V_2$.
This also shows that $G_2$ acts faithfully on $S(V)^{H_1}/v_1S(V)^{H_1}$.
Indeed if $g \in G$ with $(g-1)(S(V)^{H_1}) \subseteq v_1S(V)^{H_1}$, then $(g-1)(u_i)^{p^{f_i}}\in v_1S(V)$ and since the latter is prime ideal $(g-1)(u_i)\in v_1S(V)\cap V=Fv_1$, implying, since $\{v_1,u_2,u_3,u_4\}$ is a basis of $V$, that $g \in H_1$. Consequently by  Proposition \ref{P2-11}, using the Cohen-Macaulay and $R_2$ properties, we have that $S(V)^G/(v_1)=[S(V)^{H_1}/v_1S(V)^{H_1}]^{G_2}$.

We may assume that $f_1\le f_2 \le f_3$.
In order to exploit the $3$-dimensional result of $S(V_2)^{G_2}$ (Theorem A), $\{\bar{u}_2,\bar{u}_3,\bar{u}_4\}$ should have been a triangulating basis for the action of $G_2$ on $V_2$.
We therefore change, by Proposition \ref{P2-9}, the generators $\{x_2,x_3,x_4\}$ into $\{y_2,y_3,y_4\}$ having the properties:
\begin{enumerate}
  \item $S(V)^{H_1}=F[v_1,y_2,y_3,y_4]$, $\deg y_i=p^{f_i}$, $y_i - w_i^{p^{f_i}}\in v_1S(V)$, $i=2,3,4$, and
  \item $(g-1)(y_2) \in Fv_1^{p^{f_2}}$, $(g-1)(y_i) \in F[v_1,y_2,...,y_{i-1}]_{+}$, $i=3,4$, $\forall g \in G$.
\end{enumerate}
Hence $(g-1)(\bar{w}_2)^{p^{f_2}}=0$, $(g-1)(\bar{w}_i)^{p^{f_i}} \in (\bar{w}_2,...,\bar{w}_{i-1})S(V_2)$, $i=3,4$, $\forall g \in G_2$.
Since the latter is a prime ideal in $S(V_2)$, we have by degree consideration:
$$g(\bar{w}_2)=\bar{w}_2, \text{ }  (g-1)(\bar{w}_i)\in (\bar{w}_2,...,\bar{w}_{i-1})S(V_2)\cap V_2
= F\bar{w}_2+\cdots +F\bar{w}_{i-1}, \text { } i=3,4, \text { } \forall g\in G_2.$$
$G_2$ is generated by transvections by Lemma \ref{L2-6}.
We can use now the $3$-dimensional result (Theorem A) according to the following separation into cases, where all matrices are given with respect to the $G_2$-triangulating basis $\{\bar{w}_2,\bar{w}_3,\bar{w}_4\}$ of $V_2$.\\
$$\textbf{1: } G_2 \subset
\begin{bsmallmatrix}
 1 & 0 & F \\
 0 & 1 & F \\
 0 & 0 & 1
\end{bsmallmatrix},
\textbf{ 2: } G_2 \subset
\begin{bsmallmatrix}
 1 & F & F \\
 0 & 1 & 0 \\
 0 & 0 & 1
\end{bsmallmatrix},
\textbf{ 3: } G_2 \not\subset
\begin{bsmallmatrix}
 1 & 0 & F \\
 0 & 1 & F \\
 0 & 0 & 1
\end{bsmallmatrix},
\begin{bsmallmatrix}
 1 & F & F \\
 0 & 1 & 0 \\
 0 & 0 & 1
\end{bsmallmatrix}.$$

\textbf{Case 1:} This case is translated into $(g-1)(\bar{w}_4)$, $(g-1)(\bar{w}_3)\in F\bar{w}_2 $, $g(\bar{w}_2)=\bar{w}_2$ $\forall g \in G_2$.
Therefore by Corollary \ref{C2-8}:
$$S(V)^G/(v_1)= F[\bar{y}_2,\bar{y}_3,\bar{y}_4]^{G_2}
=F[\bar{w}_2^{p^{f_2}},\bar{w}_3^{p^{f_3}},\bar{w}_4^{p^{f_4}}] ^{G_2} =F[\bar{w}_2^{p^{f_2}},y,z].$$
Taking $\{x',y',z'\}\subset S(V)^G$, (homogenous) preimages of $\{\bar{w}_2^{p^{f_2}},y,z\}$, we get that $S(V)^G=F[v_1,x',y',z']$ as needed.

\textbf{Case 2:} This case is translated into  $g(\bar{w}_2)=\bar{w}_2$, $g(\bar{w}_3)=\bar{w}_3$, $(g-1)(\bar{w}_4)\in F\bar{w}_2+F\bar{w}_3$, $\forall g \in G_2$. Therefore by Proposition \ref{P2-6.5} $F[\bar{w}_2,\bar{w}_3,\bar{w}_4]^{G_2}=F[\bar{w}_2,\bar{w}_3,a]$, where $\deg a=|G_2|$. Therefore, using $f_2\le f_3 \le f_4$, $\{\bar{w}_2^{p^{f_2}},\bar{w}_3^{p^{f_3}},a^{p^{f_4}}\}$ is a s.o.p. in $S(V)^G/(v_1)=F[\bar{w}_2^{p^{f_2}},\bar{w}_3^{p^{f_3}},\bar{w}_4^{p^{f_4}}]^{G_2}$.
Taking $\{x',y',z'\}$ in $S(V)^G$, which are homogenous preimages of $\{\bar{w}_2^{p^{f_2}},\bar{w}_3^{p^{f_3}},a^{p^{f_4}}\}$, we get that $\{v_1,x',y',z'\}$ is a s.o.p. in $S(V)^G$ and:
$$(\deg v_1)(\deg x')(\deg y')(\deg z')= 1\cdot p^{f_2} \cdot p^{f_3} \cdot p^{f_4}(\deg a) = |H_1||G_2|=|H_1||G/H_1|=|G|.$$
Hence by Proposition \ref{P2-11.5} $S(V)^G=F[v_1,x',y',z']$.

\textbf{Case 3:} By Theorem \ref{T3-3} $S(V_2)^{G_2}=F[\bar{w}_2,\bar{w}_3,\bar{w}_4]^{G_2}=F[\bar{w}_2,a,b]$, where $a \in F[\bar{w}_2,\bar{w}_3]^{G_2}$.
Hence $(\deg a)(\deg b)=|G_2|$. The argument now follows as in the previous case, using $f_2\le f_3 \le f_4$.
We have that $\{\bar{w}_2^{p^{f_2}},a^{p^{f_3}},b^{p^{f_4}}\}$ is a s.o.p. in $F[\bar{w}_2^{p^{f_2}},\bar{w}_3^{p^{f_3}},\bar{w}_4^{p^{f_4}}]^{G_2}=S(V)^G/(v_1)$.
Taking $\{x',y',z'\}$, homogenous preimages in $S(V)^G$ of $\{\bar{w}_2^{p^{f_1}}, a^{p^{f_3}}, b^{p^{f_4}}\}$, we get that $\{v_1,x',y',z'\}$ is a s.o.p. in $S(V)^G$ with:
$$(\deg v_1)(\deg x')(\deg y')(\deg z')= 1 \cdot p^{f_2}(\deg a)(p^{f_3})(\deg b)p^{f_4}=|H_1||G_2|=|G|.$$
Hence $S(V)^G=F[v_1,x',y',z']$.


\section{\bf Acknowledgements}\label{Acknowledgements}
I thank Gregor Kemper for discussions regarding \cite{K-Loci} and Lemma \ref{L2-4}.

\end{document}